\documentclass[12pt]{article} 

\usepackage{graphicx} 
\usepackage{amsmath,amssymb}
\usepackage[all]{xy} 
\usepackage{tikz-cd}
\usetikzlibrary{arrows.meta,positioning}
\usepackage{tocloft}
\usepackage{url}       
\usepackage{setspace}
\usepackage{amsthm}
\newtheorem{theorem}{Theorem}[section]
\usepackage{graphicx}
\theoremstyle{plain}
\newtheorem*{corollary*}{Corollary}
\newtheorem*{theorem*}{Theorem}
\usepackage{hyperref}
\usepackage{amsmath,amssymb}
\newcommand{\KO}{KO}

\usepackage{slashed}  

\usepackage{amsthm}

\theoremstyle{definition}

\newtheorem*{definition*}{Definition}

\usepackage{amsthm}

\theoremstyle{remark}
\newtheorem*{remark*}{Remark}

\usepackage{amsmath,amssymb,amsthm}
\usepackage{enumitem}

\setlist[enumerate,1]{label=(\roman*)}

\newcommand{\CS}{\operatorname{CS}}

\newcommand{\Ahat}{\widehat{A}}

\DeclareMathOperator{\Hom}{Hom}

\newtheorem*{proposition*}{Proposition}

\newtheorem*{lemma*}{Lemma}

\begin{document}

\title{ Index invariants and Eta invariants determine Differential KO theory in degrees that are multiples of 
8}
\author{Tan Su }
\date{}

\maketitle

\begin{abstract}
Sullivan--Simons developed a Cheeger--Simons differential character analogue
for degree \(0 \bmod 2\) differential \(K\)-theory, giving a complete set of numerical invariants
that determine a complex vector bundle on a base manifold \(X\) with unitary connection,
up to Chern--Simons equivalence of the connection.
In this paper we develop a degree \(0 \bmod 8\) differential \(KO\) analogue. Namely,
given a real vector bundle with orthogonal connection, we construct
\(\mathbb{R}/\mathbb{Z}\)-valued \(\eta\)-invariants in the context of Atiyah--Patodi--Singer and
\(\mathbb{Z}_{2}\) Atiyah--Singer index invariants
that completely determine differential \(KO\) in degree \(0 \bmod 8\); we call this
the \emph{differential \(KO\) character}.

\end{abstract}
\section{Introduction}

Let \(M\) be a closed \(\mathrm{Spin}\) manifold mapping smoothly \(f\colon M\to X\) to a target manifold \(X\).
Suppose \(X\) carries a real vector bundle \(E\) with an orthogonal connection \(\nabla\).
Assume \(M\) is equipped with a Riemannian metric, inducing the Levi–Civita connection on \(TM\).

If \(M\) has dimension \(4k-1\), one can define an invariant given by the (reduced) eta invariant modulo \(\mathbb{Z}\) of the Spin Dirac operator on \(M\) twisted by the pullback bundle \(f^{*}E\).

If \(M\) has dimension \(8k+1\) or \(8k+2\), one can define a \(\mathbb{Z}_{2}\)-valued invariant via the Wrong Way map constructed by Atiyah-Singer ,
\[
KO^{0}(M)\longrightarrow KO^{-(8k+1)}(\mathrm{pt})\ \text{ or }\ KO^{-(8k+2)}(\mathrm{pt}),
\]
which is the index of the Spin Dirac operator twisted by the pullback bundle \(f^{*}E\).

Let $M^{n}$ be a closed spin manifold with $n=8k+1$ or $n=8k+2$, and let
$g:M\to \mathrm{pt}$ be the constant map.
The spin structure gives a $KO$–orientation of the stable normal bundle
$\nu\to M$, hence a Thom class
$u\in KO^{r}_{c}(\nu)$, where $r=\operatorname{rank}(\nu)$, and a Thom
isomorphism
\[
KO^{0}(M)\;\xrightarrow{\ \cong\ }\;KO^{r}_{c}(\nu).
\]
Choose an embedding $M\hookrightarrow S^{N}$ with normal bundle $\nu$; the
Pontrjagin–Thom collapse is a map of spaces
\[
S^{N}\longrightarrow \operatorname{Th}(\nu),
\]
which on $KO$–cohomology is implemented as extension by zero
\[
KO^{r}_{c}(\nu)\longrightarrow KO^{r}(S^{N}).
\]
Composing Thom isomorphism with extension by zero and using Bott periodicity
$KO^{r}(S^{N})\cong KO^{r-N}(\mathrm{pt})=KO^{-n}(\mathrm{pt})$ gives the
wrong–way map
\[
g_{!}:KO^{0}(M)\longrightarrow KO^{-n}(\mathrm{pt}).
\]
Since
\[
KO^{-n}(\mathrm{pt})\cong
\begin{cases}
KO^{-1}(\mathrm{pt})\cong\mathbb{Z}_{2}, & n=8k+1,\\[2pt]
KO^{-2}(\mathrm{pt})\cong\mathbb{Z}_{2}, & n=8k+2,
\end{cases}
\]
this wrong-way map defines a $\mathbb{Z}_{2}$–valued invariant, the
$KO$–index  of the spin Dirac operator on $M$.

We refer to these invariants collectively as the differential \(KO\)-character ($\widehat{KOCH}$) of \((E,\nabla)\).

One sees that the differential \(KO\)-characters of \((E,\nabla^{0})\) over \(X\) agree with those of \((E,\nabla^{1})\) whenever \(\nabla^{0}\) and \(\nabla^{1}\) are Chern--Simons (CS) equivalent. This means the odd--dimensional CS forms whose exterior derivative is the difference of the corresponding Pontryagin forms are not only closed—so that the Pontryagin character forms of \((E,\nabla^{0})\) and \((E,\nabla^{1})\) coincide—but also exact. The equality splits into two cases:

\smallskip
\noindent\emph{Case \(\dim M=4k-1\).} Use the variational formula for the spectral flow between the two Dirac operators, which is locally computable in terms of characteristic forms.

\smallskip
\noindent\emph{Case \(\dim M=8k+1,\,8k+2\).} The \(\mathbb{Z}_{2}\)--valued index is independent of the geometric data on \(E\), and hence the resulting invariants are equal.

One may also change the pair \((E,\nabla)\) by direct summing with the trivial connection on the trivial
bundle or by changing by a strict isomorphism of bundles with connection, without changing the values
of the character function associated with \((E,\nabla)\).

The goal of this paper is to show that the differential \(KO\)-character actually determines \((E,\nabla)\) up to Chern--Simons equivalence; by describing an associated hexagon that completely characterizes the properties of such differential \(KO\)-characters.

Sullivan–Simons showed that, for differential \(K\)-theory, the Hopkins–Singer axioms can be reformulated as the Grothendieck group of stable isomorphism classes of complex bundles with \emph{unitary} connections when working with \(\mathbb{R}\)-valued differential forms. With minor modifications to [SS10], an analogous statement holds for differential \(KO\)-theory in terms of \emph{real} vector bundles equipped with \emph{orthogonal} connections.

\smallskip
The term “differential \(KO\)-character’’ is inspired by Cheeger–Simons differential characters.
A Cheeger–Simons character is a function \(f\) defined on homology cycles \(M\) with values in
\(\mathbb{R}/\mathbb{Z}\), additive under disjoint union. Such an \(f\) is not a homology invariant:
if \(M=\partial W\), then
\[
f(\partial W)\;=\;\int_{W}\omega_{f}\;\;(\mathrm{mod}\ \mathbb{Z}),
\]
for an associated differential form \(\omega_{f}\). Thus \(f\) fails to descend to homology
precisely when \(\omega_{f}\) has nonzero periods.

There is a natural surjection sending a differential character \(f\) to its curvature form
\(\omega_{f}\).
The kernel consists of those characters with \(\omega_{f}=0\), which identify with
\(\operatorname{Hom}\!\bigl(H_{*}(X),\mathbb{R}/\mathbb{Z}\bigr)\), and hence with
\(H^{*}(X;\mathbb{R}/\mathbb{Z})\) by the universal coefficient theorem.

An analogue of differential characters for differential \(KO\) can be described as follows.
Replace homology cycles by spin bordism cycles enriched with a metric, and replace variation
by homology with variation by cobordism. By the Atiyah–Patodi–Singer theorem, such a cobordism
produces a real-valued form representing the Pontryagin character of a vector bundle, wedged with
the \(\widehat{A}\)-genus, integrated over the interior of the cobordism. There also arises a product
relation for differential \(KO\)-characters, analogous to  [SS18], governed by the
action of spin bordism to a point together with the \(\widehat{A}\)-genus via the Atiyah–Bott–Shapiro
orientation.

Differential \(KO\)-characters are, in general, not cobordism invariants. 
However, just as the subset of differential characters with vanishing curvature (“integral”) forms descends to homology invariants, the differential \(KO\)-characters whose associated Pontryagin forms vanish do descend to cobordism invariants. 
This is analogous to the flat index bundle phenomenon in [APS76].

The second part of this paper develops a equivalence between
family index theorems in the setting of differential $KO$ characters and
in the setting of differential $KO$ theory formulated via structured
bundles. 

For a Riemannian submersion $\pi:X\to B$ whose fibers are closed
$8k$-dimensional spin manifolds, we construct pushforwards in both theories.
In the end, the fact that these two pushforwards commute is shown to be
equivalent to the famous Bismut--Cheeger adiabatic limit theorem[BC89], which
yields a rather surprising result.

\smallskip
\noindent
(1) Via a simple isomorphism that converts the total odd form in the
Freed--Lott model of differential $KO$ to a stably trivial bundle,
we transport the Freed--Lott analytic index in [FL10] to the structured
bundle model, thereby obtaining a pushforward in the structured bundle model
of differential $KO$.

\smallskip
\noindent
(2) For the differential \(KO\) version, we define the pushforward of a
differential \(KO\)-character by evaluation on \emph{enriched spin cycles}:
for a cycle downstairs \((f:M\to B)\),
\[
(\pi_!\widehat\kappa)(M,f)\ :=\ \widehat\kappa\big(\text{pullback enriched cycle on }X\big).
\]
Given \(\pi:X\to B\) and \(f:M\to B\), form the fiber product
\[
\widetilde{M}\ :=\ M\times_{B}X\ =\ \{(m,x)\in M\times X \mid f(m)=\pi(x)\}.
\]
The pulled-back enriched cycle upstairs is the fibration
\(\widetilde M\to M\) with spin fiber \(F\) of dimension \(8k\), equipped with
the direct–sum metric. Under adiabatic scaling \(
g_{\widetilde M,t}=g_M\oplus t^{-2}g_{X/B}\), one has the limit
\[
\lim_{t\to 0}\,\bar{\eta}\!\left(D^{\,\widetilde{f}^{*}(E,\nabla)}_{\widetilde{M},t}\right).
\]

\smallskip
\noindent
We will show that the two pushforwards agree under the identification of
models, i.e. the diagram
\[
\begin{tikzcd}
\widehat{KO}(X) \arrow[d, "\check{\pi}_!^{\,\mathrm{an}}"'] \arrow[r, "\cong"] &
\widehat{KOCH}(X) \arrow[d, "\pi_!"] \\
\widehat{KO}(B) \arrow[r, "\cong"] &
\widehat{KOCH}(B)
\end{tikzcd}
\]
commutes. This yields a conceptual explanation of the Cheeger–Bismut
adiabatic limit theorem for the \(\eta\)-invariant in this setting.

The organization of the paper is as follows.

In \textbf{Section 2} we introduce the structured bundle formulation of
differential $KO$-theory in degree $0$ . In this
section we will only prove Theorem~2.6 and Proposition~2.10; everything else
follows directly, line by line, from [SS10]

In \textbf{Section 3} we show that the structured bundles in Section~2 can be
viewed as representatives of triples $(E,\nabla,\eta)$, where $\eta$ is a
total odd form in the Freed--Lott model of differential $KO$-theory, in the
special case where $\eta$ is taken to be $0$. The purpose of this section is
to transfer the analytic index of [FL10] directly to our model.

\textbf{Section 4} gives the definition of the differential $KO$ character and
explains the properties of these invariants.

\textbf{Section 5} characterizes $KO$-theory completely in terms of pairings
with Spin bordism cycles with compatible $\mathbb{Q}$ and
$\mathbb{Q}/\mathbb{Z}$-periods (inscribed in index theory) over closed
manifolds and Sullivan's  $\mathbb{Z}_{k}$-manifolds in the sense of 
and [MS74],[Fre88],[Hu23]. One of the main theorems of this section is that our
differential $KO$ character induces such $\mathbb{Q}$ and
$\mathbb{Q}/\mathbb{Z}$-valued period maps thus determining a real vector bundle up to stably equivalence. Along the way, we obtain an
interesting corollary that provides a complete characterization of when a form
is cohomologous to the Pontryagin character of some vector bundle.

\textbf{Section 6} shows that differential $KO$ characters are isomorphic to
the structured bundle model of differential $KO$ via an isomorphism of the hexagon diagram, which completely determines differential $KO$ at degree 8k.

In \textbf{Section 7} we develop family index theorems for both the
differential $KO$ character and the structured bundle model, and arrive at an
interesting discussion: the statement that these two wrong-way maps``commute''
is equivalent to the famous Cheeger--Bismut adiabatic limit theorem[BC89].

In the \textbf{final section (Section 8)} of the paper, we describe how our
results relate to classical developments from the 1970s, including Adams's
$e$-invariant on framed $(4k-1)$-manifolds, the relationship with the
Cheeger--Simons Pontryagin classes and their associated invariants on
$3$-manifolds, and cobordism invariants for flat bundles discovered by
Atiyah--Patodi--Singer.

In \textbf{Appendix A} we present a natural procedure for inducing a
Cheeger--Simons differential character to a $\mathbb{Z}_{k}$-valued period map
on cohomology. This example suggests that the Morgan--Sullivan description of
ordinary cohomology---as the group of ``liftable'' $\mathbb{Q}/\mathbb{Z}$-valued
period maps on homology with $\mathbb{Q}/\mathbb{Z}$ coefficients is already compatible with the original Cheeger--Simons character.

In \textbf{Appendix B} we briefly discuss the phenomenon that the differential
Euler character in odd dimensions can be independent of geometric data (so the
associated curvature form vanishes). This example illustrates the same
behavior encountered in our $KO$-character on enriched Spin cycles of
dimensions $8k+1$ and $8k+2$.
\section*{Acknowledgments}
The author is deeply grateful to Dennis Sullivan for many helpful conversations and for
his constant inspiration. The author also wishes to thank Jeff Cheeger for introducing the author to differential characters and for his insights on the $\eta$-invariant. The author is further indebted to Dan Freed, Daniel Grady, Jiahao Hu, Ed Miller, and Mahmoud Zeinalian for helpful conversations and valuable advice.

\section{Structured Bundle in Differential KO theory}

 From now on, \(\Lambda^{4k-1}\) denotes \(\displaystyle \bigoplus_{m\ge 0}\Omega^{4m-1}(M)\).

\[
H^{4k}(X) \;:=\; \bigoplus_{m\ge 0} H^{4m}(X),
\qquad
H^{4k-1}(X) \;:=\; \bigoplus_{m\ge 0} H^{4m-1}(X).
\]

We will continue to use this notation for the remainder of this section.

We note that most statements and arguments parallel [SS10]; many propositions will therefore be stated without proof, and interested readers are referred to [SS10] for details. In this section we prove only Theorem~2.6  and Proposition~2.10 and we later use them to construct an isomorphism between the Freed–Lott model of differential \(KO\) and the structured–bundle model of Sullivan–Simons.
\medskip

\textbf{Definition 2.1     } Let $V\to X$ be a real vector bundle with Euclidean metric $\langle\ ,\ \rangle$.
A \emph{metric connection} $\nabla$ on $V$ means $\nabla$ is compatible with
$\langle\ ,\ \rangle$ (i.e. an $O(\mathrm{rank}\,V)$--connection).
Write $R^\nabla$ for its curvature.

The \emph{Pontryagin character form} of $(V,\nabla)$ is the closed even form
\[
Ph(\nabla)\ :=\ \sum_{k\ge 0}\ \frac{(-1)^k}{(2\pi i)^{2k}(2k)!}\,\mathrm{tr}\!\big((R^\nabla)^{2k}\big)
\ \in\ \Omega^{\mathrm{4k}}_{\mathrm{cl}}(X),
\]
whose de Rham class equals the Pontryagin character $ph(V)\in H^{4\bullet}(X;\mathbb{R})$
(equivalently, $Ph(\nabla)=ch(\nabla^{\mathbb{C}})$ restricted to degrees $4k$).

Given a smooth path \(\gamma=\{\nabla^{t}\}_{t\in[0,1]}\) of metric connections,
let \(\widetilde{\nabla}\) be the induced connection on \([0,1]\times X\) and set
\begin{equation*}
cs_{Ph}(\gamma)\ :=\ \int_{0}^{1}\iota_{\partial_{t}}\,Ph(\widetilde{\nabla})\,dt
\ \in\ \Omega^{4k-1}(X).
\end{equation*}
Then
\begin{equation*}
d\,cs_{Ph}(\gamma)\ =\ Ph(\nabla^{1})-Ph(\nabla^{0}).
\end{equation*}

\paragraph{Proposition 2.2} [SS10]
If \(\alpha\) and \(\gamma\) are two paths of metric connections from \(\nabla^{0}\) to \(\nabla^{1}\), then
\[
cs_{Ph}(\alpha)\;=\;cs_{Ph}(\gamma)\quad\text{mod exact}.
\]
Consequently we may define
\begin{equation*}
CS_{Ph}(\nabla^{0},\nabla^{1})\ :=\ cs_{Ph}(\gamma)\ \ \text{in}\ \
\Omega^{4k-1}(X)/\mathrm{im}\,d,
\end{equation*}
for any such path \(\gamma\), and the previous identity becomes
\[
d\,CS_{Ph}(\nabla^{0},\nabla^{1})\ =\ Ph(\nabla^{1})-Ph(\nabla^{0}).
\]

Moreover, for any smooth map \(\psi:Y\to X\) and metric connections \(\nabla^{0},\nabla^{1}\) on \(V\),
\begin{equation*}
\psi^{*}CS_{Ph}(\nabla^{0},\nabla^{1})\;=\;CS_{Ph}(\psi^{*}\nabla^{0},\psi^{*}\nabla^{1}).
\end{equation*}

\paragraph{Definition 2.3 }
For $O(n)$ connections \(\nabla^{0},\nabla^{1}\) on \(V\), write
\[
\nabla^{0}\ \sim_{CS}\ \nabla^{1}
\quad\Longleftrightarrow\quad
CS_{Ph}(\nabla^{0},\nabla^{1})=0\ \ \text{in}\ \ \Omega^{\mathrm{4k-1}}(X)/\mathrm{im}\,d.
\]

\paragraph{Definition 2.4}
A \emph{structured (real) bundle} over \(X\) is a pair
\[
\mathcal{V}\ =\ [\,V,\{\nabla\}\,],
\]
where \(V\to X\) is a real vector bundle with metric and \(\{\nabla\}\) is an
equivalence class of metric connections under \(\sim_{CS}\).

\paragraph{Definition 2.5}
Two structured bundles \(\mathcal{V}=[V,\{\nabla\}]\) and \(\mathcal{W}=[W,\{\nabla^{W}\}]\)
are \emph{isomorphic} if there exists a bundle isometry \(\sigma:V\!\to W\)
covering \(\mathrm{id}_{X}\) such that \(\sigma^{*}\{\nabla^{W}\}=\{\nabla\}\).

\paragraph{Theorem 2.6 }
Given any \(\mathcal{V}\in \mathrm{Struct}(X)\) there is a \(\mathcal{W}\in \mathrm{Struct}(X)\) such that
\(\mathcal{V}\oplus \mathcal{W} = [n]\) for some \(n\).
Any such \(\mathcal{W}\) will be called an \emph{inverse} of \(\mathcal{V}\).

\begin{proof}
Let\[B_kO(n)\;\cong\;\mathrm{Gr}_n(R^{\,n+k})
\]
be the real Grassmannian, with tautological bundles
\[
0\longrightarrow \gamma^n \longrightarrow \varepsilon^{\,n+k}
\longrightarrow \gamma^k \longrightarrow 0,
\qquad
\gamma^n \oplus \gamma^k \;\cong\; \varepsilon^{\,n+k}.
\]
Equip the trivial bundle \(\varepsilon^{\,n+k}\) with its flat metric connection \(d\).
Let \(\nabla^n\) and \(\nabla^k\) be the metric (orthogonal–projection) connections
on \(\gamma^n\) and \(\gamma^k\) induced from \(d\). Then
\[
\nabla^n \oplus \nabla^k \;=\; d
\quad\text{on}\quad \gamma^n \oplus \gamma^k \cong \varepsilon^{\,n+k}.
\]

Recall lemma 1.16 from [SS10]
If \((E,\nabla^E)\oplus(F,\nabla^F)\cong(\varepsilon,d)\) with \(d\) flat, then
their structured classes satisfy
\[
[E,\{\nabla^E\}] + [F,\{\nabla^F\}] \;=\; 0
\quad.\] under equivalence of structured bundle.

\[
[\gamma^n,\{\nabla^n\}] + [\gamma^k,\{\nabla^k\}] \;=\; 0,
\qquad
\]

By existence of universal connection from Narasimhan--Ramanan (Applied tocompact lie group. \(O(n)\)), for sufficiently large \(k\) and any rank-\(n\)
 real bundle with O(n) connection \(\mathcal V\) over a smooth manifold \(X\), there is a smooth classifying map
\[
f: X \longrightarrow B_kO(n)
\quad\text{such that}\quad
\mathcal V \;\cong\; f^{*}[\gamma^n,\{\nabla^n\}].
\]
The inverse of \(\mathcal V\) is obtained by pulling back  $[\gamma^k,\{\nabla^k\}]$
\end{proof}

\paragraph{Definition 2.7}
Using the Grothendieck construction that turns a commutative semigroup into an abelian group, define
\[
\widehat{KO}\;=\;KO\bigl(\mathrm{Struct}(X)\bigr).
\]
Explicitly, \(\widehat{KO}(X)\) is the free abelian group on isomorphism classes of structured bundles, subject to the relation
\([\mathcal{V}]+[\mathcal{W}]=[\mathcal{V}\oplus\mathcal{W}]\).
\medskip

\paragraph{Definition 2.8}
\[
\mathrm{Struct}_{\mathrm{ST}}(X)
=\bigl\{\, [V,\{\nabla\}] \in \mathrm{Struct}(X)\ \big|\ V \text{ is stably trivial}\,\bigr\}.
\]
\[
\mathrm{Struct}_{\mathrm{SF}}(X)
=\bigl\{\, [V,\{\nabla\}] \in \mathrm{Struct}(X)\ \big|\ V \text{ is stably flat }\,\bigr\}.
\]

\medskip

For \(\mathcal{V}\in \mathrm{Struct}_{\mathrm{ST}}(X)\), let \(F\) and \(H\) be trivial bundles such that
\(V\oplus F = H\), and let \(\nabla^{F},\nabla^{H}\) be flat connections on \(F\) and \(H\).
We define
\[
\widehat{CS}:\ \mathrm{Struct}_{\mathrm{ST}}(X)\longrightarrow
\Lambda^{\mathrm{4k-1}}/\Lambda_{O}
\]
by
\[
\widehat{CS}(\mathcal{V})
=\ CS\!\bigl(\nabla^{H},\,\nabla \oplus \nabla^{F}\bigr)
\ \ \mathrm{mod}\ \ \Lambda_{O}/\Lambda^{\mathrm{4k-1}}_{\mathrm{exact}}\, .
\]

we note 
$\wedge_{O} = \{g^{*}(\Theta)\} + \wedge^{4k-1}_{exact}$, where
$\Theta$ represents the
universal transgression  of the pontryagin character and $g: X \rightarrow
BO$ runs through all smooth maps classifying the given(stabilized) Real Bundle.

\paragraph{Proposition 2.9 }
\;\(\ker(\widehat{CS})=\mathrm{Struct}_{\mathrm{SF}}(X).\)

This proposition follows identically from Prop 2.5 [SS10]

\paragraph{Proposition 2.10 }\;
\(\operatorname{Im}(\widehat{CS})=\Lambda^{\mathrm{4k-1}}(X)\big/\Lambda_{O}(X).\)

\begin{proof}

We first show \(\Lambda^{3}(\mathbb{R}^{n})\big/\Lambda_{O}(\mathbb{R}^{n})
\subseteq \operatorname{Im}(\widehat{CS})\).
Let \(\eta=f(x)\,dx_i\wedge dx_j\wedge dx_k\in\Omega^3(\mathbb{R}^n)\) with \(i<j<k\).
On the trivial rank-\(2\) real bundle \(E\cong\mathbb{R}^{2}\), pick a \(1\)-form \(w\) and set
\[
A=\begin{pmatrix}0 & w\\[2pt] -w & 0\end{pmatrix}\in \Omega^1(\mathbb{R}^n;\mathfrak{so}(2))
\subset \Omega^1(\mathbb{R}^n;\mathfrak{so}(n)).
\]
Then \(A\wedge A=0\) and the curvature is
\(F=dA=\begin{pmatrix}0 & dw\\ -dw & 0\end{pmatrix}\).
Along the straight path \(\nabla^t=tA\), the Pontryagin–Chern–Simons form for \(p_1\) is
\[
\widehat{CS}(\nabla^t)=c_1\,w\wedge dw \qquad (\mathrm{mod}\ d\Omega^{2}),
\]
for a fixed nonzero constant \(c_1\).
Note
\[
\bigl(\widehat{CS}(\nabla ^{t})\bigr)_{4k-1}
\;=\;
\frac{1}{k\,(2k-1)!\,(2\pi)^{2k}}\;
\omega \wedge (d\omega)^{\,2k-1}.
\]

Choose
\[
w := x_i\,dx_j + f\,dx_k .
\]
Then \(dw=dx_i\wedge dx_j + df\wedge dx_k\) and
\begin{align*}
w\wedge dw
&=(x_i\,dx_j+f\,dx_k)\wedge(dx_i\wedge dx_j+df\wedge dx_k)\\
&= \underbrace{x_i\,dx_j\wedge dx_i\wedge dx_j}_{0}
   \;+\; x_i\,dx_j\wedge df\wedge dx_k
   \;+\; f\,dx_k\wedge dx_i\wedge dx_j
   \;+\; \underbrace{f\,dx_k\wedge df\wedge dx_k}_{0}.
\end{align*}
Using
\[
d\!\big(x_i f\,dx_j\wedge dx_k\big)
= f\,dx_i\wedge dx_j\wedge dx_k + x_i\,df\wedge dx_j\wedge dx_k,
\]
we get
\(x_i\,dx_j\wedge df\wedge dx_k \equiv f\,dx_i\wedge dx_j\wedge dx_k\) modulo exact forms.
Hence
\[
w\wedge dw \equiv 2\,f\,dx_i\wedge dx_j\wedge dx_k \quad (\mathrm{mod}\ d\Omega^{2}),
\]
and therefore
\[
\widehat{CS}(\nabla^t)=c_1\,w\wedge dw 
\ \equiv\ (2c_1)\,f\,dx_i\wedge dx_j\wedge dx_k 
\quad (\mathrm{mod}\ exact).
\]
Since \(2c_1\neq 0\), this produces the desired basic \(3\)-form modulo exact forms,
hence modulo \(\Lambda_O\).

Fix distinct indices \(i_1<\cdots<i_{4k-1}\) and set
\[
w \;=\; \sum_{r=1}^{2k-1} x_{i_{2r-1}}\,dx_{i_{2r}} \;+\; f\,dx_{i_{4k-1}},\]
Write \(\gamma:=\sum_{r=1}^{2k-1}dx_{i_{2r-1}}\wedge dx_{i_{2r}}\), so that
\(dw=\gamma+df\wedge dx_{i_{4k-1}}\).
A standard wedge computation yields
\[
w\wedge (dw)^{2k-1} \;\equiv\; (2k-1)!\,f\,dx_{i_1}\wedge\cdots\wedge dx_{i_{4k-1}}
\qquad (\mathrm{mod}\ exact).
\]

Under mod \(\Lambda_{O}(\mathbb{R}^{n})\)
\smallskip
\[
\widehat{CS}\!\big(\nabla^t\\(normalized))
\;=\;
f\,dx_{i_1}\wedge\cdots\wedge dx_{i_{4k-1}}
\;+\;\theta,
\qquad 
\theta\ \in\ \sum_{m=1}^{k-1}\Lambda^{4m-1}(\mathbb{R}^{n}),
\]

By the induction hypothesis there exists a structured stably trivial real bundle
\(\mathcal V\) with \(\widehat{CS}(\mathcal V)=\theta\); hence
\(\widehat{CS}(\nabla^t\oplus \mathcal V^{-1})
= f\,dx_{i_1}\wedge\cdots\wedge dx_{i_{4k-1}}\).

Any element of \(\Lambda^{4k-1}(\mathbb{R}^{n})\) is a sum of the above  terms. 

Pulling back under this embedding, the rest of the proof follows from proposition 2.6[SS10]

\end{proof}

\paragraph{Theorem 2.11}
There is an isomorphism
\  
\[
\widehat{CS}_{\!Ph}:\ \mathrm{Struct}_{\mathrm{ST}}(X)/\mathrm{Struct}_{\mathrm{SF}}(X)\ \cong\ \Lambda^{\mathrm{4k-1}}(X)/\Lambda_{O}(X) 
\]

 \begin{proof}
     Obtained by combining proposition 2.9 and Proposition 2.10.
 \end{proof}

\
We note that the constructions and proofs of the hexagon for structured bundles in [SS10] carry over verbatim to the real case; the reader may verify this if desired.

\vspace{.5cm}
\begin{center}
\setlength{\unitlength}{0.5cm}
\begin{picture}(24,16)(2.5,0)
\thicklines
\put(5,1){$0$}
\put(20.5,1){$0$}

\put(6,2){\vector(1,1){1.5}}
\put(18,3.5){\vector(1,-1){1.5}}

\put(8,4.5){$\wedge^{4k-1}/\wedge_{O}$}
\put(12,4.5){\vector(1,0){2.5}}
\put(16.5,4.5){$\wedge_{BO}$}
\put(13,5){\small{$d$}}
 
\put(6.5,7.5){\vector(1,-1){1.5}}
\put(10.5,7){\small{$i_{2}$}}
\put(10.5,6){\vector(1,1){1.5}}
\put(14.5,7.5){\vector(1,-1){1.5}}
\put(15.5,7){\small{$ph$}}
\put(18.5,6){\vector(1,1){1.5}}

\put(0,8){1)}
\put(3,8){$H^{4k-1}(\mathbb{R})$}
\put(12.75,8){$\hat{KO}$}
\put(20,8){$H^{4k}(\mathbb{R})$}

\put(6,9.5){\vector(1,1){1.5}}
\put(10.5,11){\vector(1,-1){1.5}}
\put(11.5,10.5){\small{$i_{1}$}}
\put(14,10.5){\small{$\delta_{2}$}}
\put(14.0,9.5){\vector(1,1){1.5}}
\put(18,11){\vector(1,-1){1.5}}

\put(7,12){$KO^{-1}(\mathbb{R}/\mathbb{Z})$}
\put(12,12){\vector(1,0){2.5}}
\put(15.5,12){$KO^{}(\mathbb{Z})$}

\put(5.5,14.5){\vector(1,-1){1.5}}
\put(18,13){\vector(1,1){1.5}}

\put(4.5,15){$0$}
\put(20,15){$0$}
\end{picture}
\end{center}

This hexagon is in the same spirit as the Simons--Sullivan construction in the complex case, with the Chern character replaced by the Pontryagin character and \(BU\) replaced by \(BO\).

We briefly recall the meaning of the symbols.  
The group \(\wedge_{BO}\) is the ring of total closed \(4k\)-forms on \(X\)
that are cohomologous to the Pontryagin character classes of real vector bundles over \(X\).
Moreover,
\[
\wedge_{O} \;=\; \{\,g^{*}(\Theta)\,\} \;+\; \wedge^{4k-1}_{\mathrm{exact}},
\]
where \(\Theta\) represents the universal transgression of the Pontryagin character and
\(g\colon X \to O\) ranges over all smooth maps classifying the given (stabilized) real bundle.

The map \(\mathrm{ph}\) sends an element of \(\widehat{KO}(X)\) to its Pontryagin character form.
Finally, the flat theory \(KO^{-1}(\mathbb{R}/\mathbb{Z})\) is the homotopy fibre obtained by
smashing the Moore spectrum with the \(KO\) spectrum.

We also note that the upper row is given by the Bockstein sequence associated to the
coefficient sequence
\[
0 \longrightarrow \mathbb{Z} \longrightarrow \mathbb{R} \longrightarrow \mathbb{R}/\mathbb{Z} \longrightarrow 0,
\]
which yields
\[
KO^{-1}(X;\mathbb{R}/\mathbb{Z}) \;\longrightarrow\; KO^{0}(X;\mathbb{Z}) \;\longrightarrow\; KO^{0}(X;\mathbb{R}),
\]
where \(KO^{0}(X;\mathbb{R})\) is, by the Chern--Dold/Pontryagin character, isomorphic to
\(H^{4k}(X;\mathbb{R})\).

\section{Freed-Lott's Differential KO }

In [Fre25], it is pointed out that Differential KO theory can be written as the homotopy pull back square in such form. 
\[
\begin{tikzcd}[column sep=large,row sep=large]
\widetilde{KO}^{\,q}(M) \arrow[r, "\text{curvature}"] \arrow[d, "\pi_{0}"']
& \Omega\bigl(M;\,\mathbb{R}[v,v^{-1}]\bigr)^{q}_{\mathrm{closed}} \arrow[d, "\text{de Rham}"] \\
KO^{q}(M) \arrow[r]
& H\bigl(M;\,\mathbb{R}[v,v^{-1}]\bigr)^{q}
\end{tikzcd}
\qquad(\deg v = 4).
\]

In these diagrams \(\mathbb{R}[u,u^{-1}]\cong K^{\bullet}(\mathrm{pt})\otimes\mathbb{R}\) and
\(\mathbb{R}[v,v^{-1}]\cong KO^{\bullet}(\mathrm{pt})\otimes\mathbb{R}\).

The geometric model represents a differential class by a triple \((E,\nabla,\eta)\),
in which \(E\to M\) is a orthogonal vector bundle with compatible covariant derivative,
and \(\eta\) is a differential form of total degree \(q-1\)  one should see [FL10] for details. Here \(\eta \in \Omega\bigl(M;\,\mathbb{R}[v,v^{-1}]\bigr)^{\,q-1}\).
When \(q=0\), the component \(\Omega\bigl(M;\,\mathbb{R}[v,v^{-1}]\bigr)^{\,q-1}\)
=is the direct sum of all degree \(4k-1\) forms,
\[
\Omega\bigl(M;\,\mathbb{R}[v,v^{-1}]\bigr)^{-1}
=\bigoplus_{k\ge 0}\Omega^{4k-1}(M).
\]

The goal of this section is to show the Structured bundle model in previous section is realized as Freed Lott's Model with a 0 representative on the total form  $\eta$

We quickly review the Differential KO group under Freed Lott's construction

\[
\begin{tikzcd}[column sep=large]
0 \arrow[r] & E_{1} \arrow[r, "i"] & E_{2} \arrow[r, shift left=0.35ex, "j"] & E_{3} \arrow[l, shift left=0.35ex, "s"] \arrow[r] & 0
\end{tikzcd}
\tag{1}
\]
is a split short exact sequence of real vector bundles with connections \(\nabla_i\) on \(E_i\to X\) for \(i=1,2,3\), we define the relative Chern–Simons transgression form
\(\mathrm{CS}(\nabla_{1},\nabla_{2},\nabla_{3})\in \Omega^{\mathrm{4k-1}}(X)/\operatorname{Im}(d)\) by
\[
\mathrm{CS}(\nabla_{1},\nabla_{2},\nabla_{3})
:= \mathrm{CS}\bigl((i\oplus s)^{*}\nabla_{2},\, \nabla_{1}\oplus \nabla_{3}\bigr),
\]
noting that \(i\oplus s: E_{1}\oplus E_{3}\to E_{2}\) is a vector bundle isomorphism.

\medskip

\paragraph{Definition 3.1}

The Freed–Lott differential \(KO\)-group \(\widehat{KO}_{\mathrm{FL}}(X)\) is the abelian group with the following generators and relation:
a generator of \(\widehat{KO}_{\mathrm{FL}}(X)\) is a quadruple
\(\mathcal{E}=(E,\nabla,\phi)\), where \((E,\nabla)\) is as before and
\(\phi\in \Omega^{\mathrm{4k-1}}(X)/\operatorname{Im}(d)\).
The only relation is \(\mathcal{E}_{2}=\mathcal{E}_{1}+\mathcal{E}_{3}\) if and only if there exists a short exact
sequence of orthognal real vector bundles \((1)\) and
\[
\phi_{2}=\phi_{1}+\phi_{3}-\mathrm{CS}(\nabla_{1},\nabla_{2},\nabla_{3}).
\]

For \(\mathcal{E}_{1},\mathcal{E}_{2}\in \widehat{K}_{\mathrm{FL}}(X)\), the addition
\[
\mathcal{E}_{1}+\mathcal{E}_{2}
:=\bigl(E_{1}\oplus E_{2},\, \nabla^{E_{1}}\oplus \nabla^{E_{2}},\, \phi_{1}+\phi_{2}\bigr)
\]
is well defined. Note that \(\mathcal{E}_{1}=\mathcal{E}_{2}\) if and only if there exists
\((F,\nabla^{F},\phi^{F})\in \widehat{KO}_{\mathrm{FL}}(X)\) such that
\begin{enumerate}\itemsep0.2em
\item \(E_{1}\oplus F \cong E_{2}\oplus F\), and
\item \(\phi_{1}-\phi_{2}=\mathrm{CS}\bigl(\nabla^{E_{2}}\oplus \nabla^{F},\, \nabla^{E_{1}}\oplus \nabla^{F}\bigr).\)
\end{enumerate}\

We now show the structured bundle model $\widehat{KO}_{}(X)$ admits as a representative of Freed-Lott's Differential KO model. We will later use this construction to construct an analytic pushforward in Differential KO theory.

\noindent\textbf{Proposition 3.2}

$\mathbb{R}^n,[\nabla^{d}]$ denotes trivial bundle with trivial connection.
$\widehat{CS}^{-1}$ denotes the isomorphism in theorem 2.9.

Let \(X\) be a compact manifold. The maps
\[
f:\ \widehat{KO}_{}(X)\longrightarrow \widehat{KO}_{\mathrm{FL}}(X),
\qquad
g:\ \widehat{KO}_{\mathrm{FL}}(X)\longrightarrow \widehat{KO}_{}(X)
\]
are defined by
\[
f\!\left([E,[\nabla^{E}]]-[\mathbb{R}^n,[\nabla^{d}]]\right)
  \;=\; (E,\nabla^{E},0)\;-\;(R^n,d,0),
\]
\[
g(E,\nabla^{E},\phi)
  \;=\; [E,[\nabla^{E}]]\;+\;[\widehat{CS}^{-1}(\phi)
]
      \;-\;[\mathbb{R}^{\dim \widehat{CS}^{-1}(\phi)
},[d]].
\]

is ring isomorphism, f and g are inverse to each other.
\begin{proof}

It is straightforward to check f and g are well defined. 

By Theorem~2.6, every element of \(\widehat{KO}(X)\) can be written in the form
\(V - [n]\), where \(V\) is a structured bundle and
\([n] = [\mathbb{R}^{n},[d]]\) denotes the trivial rank-\(n\) bundle equipped
with the trivial connection.

\[
\begin{aligned}
f\!\left([E,[\nabla^{E}]]-[\mathbb{R}^{n},[d]]\right)
&=(E,\nabla^{E},0)-(\mathbb{R}^{n},d,0),\\
g\!\left((E,\nabla^{E},0)-(\mathbb{R}^{n},d,0)\right)
&=\Big([E,[\nabla^{E}]]+[\widehat{CS}^{-1}(0)]-[\mathbb{R}^{\dim\widehat{CS}^{-1}(0)},[d]]\Big)\\
&\quad-\Big([\mathbb{R}^{n},[d]]+[\widehat{CS}^{-1}(0)]-[\mathbb{R}^{\dim\widehat{CS}^{-1}(0)},[d]]\Big)\\
&=[E,[\nabla^{E}]]-[\mathbb{R}^{n},[d]],
\end{aligned}
\]
since \(\widehat{CS}^{-1}(0)\) is the trivial class (rank \(0\)). Hence \(g\circ f=\mathrm{id}_{\widehat{KO}(X)}\).

Compute  \(f\circ g\)
Let \((E,\nabla^{E},\phi)\in \widehat{KO}_{\mathrm{FL}}(X)\), and choose
\(\mathcal{V}:=\widehat{CS}^{-1}(\phi)=(V,[\nabla^{V}])\in \mathrm{Struct}_{\mathrm{ST}}(X)\) with
\(m:=\dim V\).
By definition,
\[
(f\circ g)(E,\nabla^{E},\phi)
=(E,\nabla^{E},0)\;+\;(V,\nabla^{V},0)\;-\;(R^{m},d,0).
\]
it is equivalent to show
\[
(E,\nabla^{E},\phi)+(R^{m},d,0)
=\ (E,\nabla^{E},0)\;+\;(V,\nabla^{V},0),
\]
i.e.
\begin{equation}\label{eq:eq-classes}
(E\oplus R^{m},\,\nabla^{E}\!\oplus d,\,\phi)
=\ (E\oplus V,\,\nabla^{E}\!\oplus\nabla^{V},\,0)
\quad\text{in }\widehat{KO}_{\mathrm{FL}}(X).
\end{equation}

Since \(\mathcal{V}\) is stably trivial, there exist trivial structured bundles
\(\mathcal{F}=(F,[d^{F}])\) and \(\mathcal{H}=(H,[d^{H}])\) with
\(V\oplus F \cong H\) and
\[
\phi=\CS(d^{H},\,\nabla^{V}\!\oplus d^{F})\in \Omega^{4k-1}(X)/Im(d).
\]
using
additivity and naturality of the transgression, we get
\[
\CS(\nabla^{E}\!\oplus\nabla^{V},\,\nabla^{E}\!\oplus d)
=\CS(\nabla^{E}\!\oplus d^{H},\,\nabla^{E}\!\oplus\nabla^{V}\!\oplus d^{F})
=\phi.
\]
Thus \eqref{eq:eq-classes} holds, and consequently
\[
(f\circ g)(E,\nabla^{E},\phi)=(E,\nabla^{E},\phi).
\]

\end{proof}

\section{Differential KO character}

\textbf{Definition 4.1 }
Let $X$ be a compact manifold. An \emph{enriched spin $n$–cycle over $X$} is a triple
\[
(M^{n},\,f,\,g^{M}),
\]
where
\begin{itemize}
  \item $M^{n}$ is a closed smooth spin manifold,
  \item $f\colon M\to X$ is a continuous map, and
  \item $g^{M}$ is a Riemannian metric on $M$ (we denote by $\nabla^{TM}$ its
        Levi–Civita connection).
\end{itemize}

\textbf{Definition 4.2 }
An enriched spin cycle $(M^{n},f,g^{M})$ is an \emph{enriched boundary}
if there exists a compact spin $(n\!+\!1)$–manifold $W$ with boundary
$\partial W=M$, a map $F\colon W\to X$ with $F|_{M}=f$, and a Riemannian metric
$g^{W}$ on $W$ such that:
\begin{enumerate}
  \item $g^{W}$ is of \emph{product form} on a collar neighborhood
        $[0,\varepsilon)\times M\subset W$, i.e.\ $g^{W}=dt^{2}\oplus g^{M}$, and
  \item the Levi–Civita connection $\nabla^{TW}$ of $g^{W}$ restricts along
        the boundary to the Levi–Civita connection $\nabla^{TM}$ of $g^{M}$.
\end{enumerate}
In this situation we say that $(M,f,g^{M})$ is an enriched boundary
(of $(W,F,g^{W})$).

\noindent
Under disjoint union, enriched spin cycles (resp.\
enriched boundaries) over $X$ form commutative monoids, which we denote
by $EC^{\mathrm{Spin}}(X)$ and $EB^{\mathrm{Spin}}(X)$, respectively.
\medskip

\emph{Remarks.}
(1) The choice of metric (and hence its Levi–Civita connection) on the spin cycle
is part of the enrichment.

(2) The enriched spin bordism of a point acts on the enriched spin cycles of \(X\)
via the external product. More explicitly, if
\((V,g^{V})\) is an enriched spin bordism cycle over a point and
\((M,f,g^{M})\) is an enriched spin cycle over \(X\), then their natural product is
\((V\times M,\,f\circ \mathrm{pr}_{M},\,g^{V}\oplus g^{M})\) over
\(X\times \{\mathrm{pt}\}\cong X\), where \(g^{V}\oplus g^{M}\) is the natural
product (direct sum) metric on \(V\times M\).

We now give the definition of the differential \(KO\) character associated to
\((E,\nabla)\).

We define 
\[
\hat{\kappa}_{(E,\nabla^{0})}(M,g) \;:=\;
\begin{cases}
\displaystyle \bar{\eta}\!\left(D_{M}\otimes f^{*}(E,\nabla^{0})\right)\in\mathbb{R}/\mathbb{Z},
& \text{if }\dim M=8k+7,\\[8pt]
\displaystyle \tfrac{1}{2}\,\bar{\eta}\!\left(D_{M}\otimes f^{*}(E,\nabla^{0})\right)\in\mathbb{R}/\mathbb{Z},
& \text{if }\dim M=8k+3,
\end{cases}
\]
where \(\bar{\eta}(D):=\frac{\eta(D)+h(D)}{2}\ (\mathrm{mod}\ \mathbb{Z})\).

For \(\dim M=8k+1\) or \(8k+2\), set
\[
\hat{\kappa}_{(E,\nabla^{0})}(M,g)
\;:=\; \pi_{!}\!\left(f^{*}E,\nabla^{0}\right)
\in
\begin{cases}
KO^{-(8k+1)}(\mathrm{pt}),& \text{if }\dim M=8k+1,\\
KO^{-(8k+2)}(\mathrm{pt}),& \text{if }\dim M=8k+2,
\end{cases}
\]
where the Atiyah--Singer wrong-way map to a point is
\[
\pi_{!}:\ KO^{0}(M)\longrightarrow KO^{-n}(\mathrm{pt}),\qquad n=\dim M.
\]

\paragraph{Theorem 4.4}
Properties of the differential \(KO\)-character \(f\)

\begin{enumerate}[label=\textbf{(\arabic*)}, leftmargin=*, align=left]
\item \textbf{Change by Cobordism.}
Let $S$ be an enriched spin $(k+1)$ manifold with boundary with $\partial S\in EB^{k}(X)$ and
Levi–Civita connection $\nabla^{TS}$. There exists a closed differential form $\omega_{f}$ on $X$ such that, after pullback to $S$,
\[
f(\partial S)\ \equiv\
\begin{cases}
\displaystyle \int_{S}\,\omega_{f}\wedge \widehat A(\nabla^{TS})
&\text{if }\dim S\equiv 0 \pmod 8,\\[1.2ex]
\displaystyle \tfrac{1}{2}\int_{S}\,\omega_{f}\wedge \widehat A(\nabla^{TS})
&\text{if }\dim S\equiv 4 \pmod 8,
\end{cases}
\quad (\mathrm{mod}\ \mathbb{Z}) .
\]
Here $\widehat A(\nabla^{TS})$ is the $\widehat A$–form of the Levi–Civita connection on $TS$.

\item \textbf{Product Relationship.}
For any enriched spin $8$–cycle $Q\in EC^{8m}(\mathrm{pt})$ and any
$M\in EC^{k}(X)$,
\[
f(Q\times M)\;=\;\widehat A(Q)\,f(M),
\]
where $\widehat A(Q)$ denotes the $\widehat A$–genus of $Q$.
\end{enumerate}

\paragraph{Remark }
It follows that the closed form $\omega_{f}$ of property has integral periods in the following sense:
for every closed 4n \emph{spin} cycle  dimensional $W \xrightarrow{\,F\,} X$,
\[
\int_{W} F^{*}\omega_{f}\wedge \widehat A(\nabla^{TW})
\in
\begin{cases}
\mathbb{Z}, & \text{if }\dim W=8K,\\[0.4ex]
2\mathbb{Z}, & \text{if }\dim W=8K+4~
\end{cases}
\]
(we show such forms are equivalent to pontryagin forms in theorem 5.4 and denote them as by $\wedge_{integrality}$ in section 6)

Before the proof of 4.4 we need one lemma.

\paragraph{Lemma 4.5}Let \(M^{d}\) be a closed Spin manifold with \(d\in\{8K+3,\,8K+7\}\), and let \(E\to M\) be a real vector bundle. Then there exists an integer \(n\ge 1\) such that \(n\) copies of \(M\) bound a Spin manifold \(W^{d+1}\) with
\[
\partial W \;=\; \bigsqcup_{1}^{\,n} M,
\]
and there exists a real vector bundle \(E_{W}\to W\) whose restriction to each boundary component is \(E\).

\begin{proof}
This follows from the fact that
\(\Omega^{\mathrm{Spin}}_{8K+3}\!\bigl(BO(n)\bigr)\) and
\(\Omega^{\mathrm{Spin}}_{8K+7}\!\bigl(BO(n)\bigr)\)
are finite torsion groups with no free part, as seen by a direct calculation using the Atiyah--Hirzebruch spectral sequence.
\end{proof}

Now we continue the proof of theorem 4.4.

\begin{proof} [Proof of Theorem 4.4]
To verify that this pairing defines a \(KO\)-character, we check:
\begin{enumerate}

\item 
Fix \((E,\nabla)\) over \(X\). Let \(f\colon M^{4k-1}\to X\) be a cycle that is null-cobordant, so there exists \(F\colon W^{4k}\to X\) with \(\partial W=M\) extending \(f\). Then by Atiyah-Patodi-Singer
\[
\hat{\kappa}_{(E,\nabla)}(M)
\;\equiv\;
\int_{W}\widehat{A}(TW)\wedge \operatorname{ph}\!\bigl(F^{*}E,F^{*}\nabla\bigr)
\pmod{\mathbb{Z}}.
\]

The differential form associated with \(\hat{\kappa}(E,\nabla)\) is the Pontryagin character form, and it satisfies the cobordism invariance criterion.

\item 

Let \(Q^{8k}\) be a closed Spin manifold and let \(M\) be an enriched Spin cycle. We want to show, under the action of \(Q\) from external product ,
\[
\hat{\kappa}_{(E,\nabla)}(Q\cdot M)
=
\widehat{A}(TQ)[Q]\;\cdot\;\hat{\kappa}_{(E,\nabla)}(M)\;\in\;\mathbb{R}/\mathbb{Z}.
\]

By lemma 4.5

We can fill \(k\) copies of \(M\) by a spin manifold \(W\) and extend \(\bigl(f^{*}E\bigr)\!\mid_{M}\) to a bundle \((f^{*}E)_{W}\to W\). Similarly,  extend \(\bigl(f^{*}\nabla^{}\bigr)\!\mid_{M}\) to a connection \(\bigl(f^{*}\nabla^{}\bigr)\!\mid_{W}\). From now, We denote  \((f^{*}E)_{W}\) as \(E_{W}\). 
Assume there is a product neighborhood near \(\partial W\) where these extended connections are product-like. On \(W\) there are two even-degree differential forms: \(\operatorname{ph}\!\bigl(E_{W},\nabla^{E_{W}}\bigr)\) and \(\widehat{A}(W)\).

Equip \(Q\times M\) with the direct-sum (product) metric.  Consider \(\widetilde{W}:=Q\times W\), so \(\partial\widetilde{W}=k(Q\times M)\). Extend \(E|_{M}\) to \(E_{W}\to W\) and use the pulled-back bundle \(\operatorname{pr}_{W}^{*}E_{W}\to Q\times W\); take product-like extensions of the connections near \(\partial\widetilde{W}\).

By Atiyah–Patodi–Singer,
\[
\int_{Q\times W}\widehat{A}\!\bigl(T(Q\times W)\bigr)\wedge
\operatorname{ph}\!\bigl(\operatorname{pr}_{W}^{*}E_{W},\nabla^{\operatorname{pr}_{W}^{*}E_{W}}\bigr)
\;\equiv\;
k\,\bar{\eta}\!\left(D_{Q\times M}\otimes \operatorname{pr}_{M}^{*}\bigl(f^{*}(E,\nabla)\bigr)\right)
\pmod{\mathbb{Z}}.
\]
Using multiplicativity of \(\widehat{A}\) for direct-sum connections,
\[
\widehat{A}\!\bigl(T(Q\times W)\bigr)
=\widehat{A}(TQ)\wedge \widehat{A}(TW),
\]
and the fact that \(\operatorname{ph}\bigl(\operatorname{pr}_{W}^{*}E_{W}\bigr)=\operatorname{pr}_{W}^{*}\operatorname{ph}(E_{W})\), we obtain
\[
\int_{Q\times W}\widehat{A}(TQ)\wedge \widehat{A}(TW)\wedge \operatorname{pr}_{W}^{*}\operatorname{ph}(E_{W})
=\widehat{A}(TQ)[Q]\cdot \int_{W}\widehat{A}(TW)\wedge \operatorname{ph}(E_{W}).
\]
Applying the APS relation for \(W\) with \(\partial W=kM\),
\[
k\,\bar{\eta}\!\left(D_{M}\otimes f^{*}(E,\nabla)\right)
\;\equiv\;
\int_{W}\widehat{A}(TW)\wedge \operatorname{ph}(E_{W})
\pmod{\mathbb{Z}},
\]
we conclude
\[
k\,\bar{\eta}\!\left(D_{Q\times M}\otimes \operatorname{pr}_{M}^{*}\bigl(f^{*}(E,\nabla)\bigr)\right)
\;\equiv\;
\widehat{A}(TQ)[Q]\cdot
k\,\bar{\eta}\!\left(D_{M}\otimes f^{*}(E,\nabla)\right)
\pmod{\mathbb{Z}}.\]
Hence the \(KO\)-character is multiplicative under the Spin cobordism action:
\[
\hat{\kappa}_{(E,\nabla)}(Q\times M)
=
\widehat{A}(TQ)[Q]\cdot \hat{\kappa}_{(E,\nabla)}(M)\in \mathbb{R}/\mathbb{Z}.
\]

In dimension n= \(8k+1\) and \(8k+2\)
By the multiplication axiom in Atiyah-Singer [1], the multiplicative property of differential KO character follows from this commutative diagram. 
\[
\begin{tikzcd}[column sep=large,row sep=large]
KO(Q\times M) \arrow[d, "\pi^{Q}_{!}"] \arrow[dr, "\pi^{Q\times M}_{!}"] & \\
KO(M) \arrow[r, "\pi^{M}_{!}"] & KO^{-n}(\mathrm{pt}) \;\cong\; \mathbb{Z}_{2}
\end{tikzcd}
\qquad\text{so that}\qquad
\pi^{Q\times M}_{!}\;=\;\pi^{M}_{!}\circ \pi^{Q}_{!}.
\]

\end{enumerate}

\end{proof}

Let \(g^{0}_{M},g^{1}_{M}\) be Riemannian metrics on \(M\) with Levi–Civita connections
\(\nabla^{TM}_{0},\nabla^{TM}_{1}\). Let \(\nabla^{E,0},\nabla^{E,1}\) be two metric connections on \(E\) over \(X\), and consider their pullbacks \(f^{*}\nabla^{E,0},f^{*}\nabla^{E,1}\) on \(f^{*}E\to M\).
Let \(D_{M,g_{i}}\otimes (f^{*}E,f^{*}\nabla^{E,i})\) denote the Dirac operator with respect to the corresponding geometric data.

\[
\hat{\kappa}_{(E,\nabla^{E,i})}(M,g_{i})
:=\bar{\eta}\!\left(D_{M,g_{i}}\otimes (f^{*}E,f^{*}\nabla^{E,i})\right)\in\mathbb{R}/\mathbb{Z}
\qquad (i=0,1).
\]

\noindent\textbf{Proposition 4.6.}
For two sets of geometric data \(g_{0},g_{1}\) on \(M\) and \(\nabla^{TM}_{0},\nabla^{TM}_{1}\) on \(TM\), and \(\nabla^{E,0},\nabla^{E,1}\) on \(E\to X\), one has
\[
\hat{\kappa}_{(E,\nabla^{E,1})}(M,g_{1})
-\hat{\kappa}_{(E,\nabla^{E,0})}(M,g_{0})
\]
\[
\equiv\ \varepsilon(M)\Bigg[
\int_{M}\widetilde{\widehat{A}}\!\bigl(TM;\nabla^{TM}_{0},\nabla^{TM}_{1}\bigr)\wedge
\operatorname{ph}\!\bigl(f^{*}E,f^{*}\nabla^{E,0}\bigr)
\;+\;
\int_{M}\widehat{A}\!\bigl(TM,\nabla^{TM}_{1}\bigr)\wedge
\widetilde{\operatorname{ph}}\!\bigl(f^{*}E;f^{*}\nabla^{E,0},f^{*}\nabla^{E,1}\bigr)
\Bigg]
\]
\[
(\mathrm{mod}\ \mathbb{Z}) .
\]
Here \(\varepsilon(M)\) depends only on \(\dim M\):
\[
\varepsilon(M)=
\begin{cases}
0, & \text{if }\dim M=1,2,\\[2pt]
\dfrac{1}{2}, & \text{if }\dim M=3,\\[6pt]
1, & \text{if }\dim M=7.
\end{cases}
\]

Here \(\widetilde{\widehat{A}}\) and \(\widetilde{\operatorname{ph}}\) denote the Chern–Simons forms associated to \(\widehat{A}\) and \(\operatorname{ph}\), respectively.

\begin{proof}
This is a direct consequence of the variational formula for the (reduced) eta invariant . Upon passing modulo \(\mathbb{Z}\), the spectral flow contribution of the Dirac operator vanishes, yielding the stated identity. in case of dim M=8K+1 and dim M=8K+2. The associated invariants yields value in \(\mathbb{Z}_2\) and don't depend on the geometric data. 
\end{proof}

\section{Invariants of Real Bundle}
In this section we first review the construction of Anderson Duality. We arrive at a description of \(KO^{0}(X)\) in terms of a commutative square of “liftable homomorphisms” from \(KO\)-homology with \(\mathbb{Q}/\mathbb{Z}\) coefficients to \(\mathbb{Q}/\mathbb{Z}\). Geometrically, this means compatible \(\mathbb{Q}\) and \(\mathbb{Q}/\mathbb{Z}\)-periods of Real \(K\)-theory, expressed via the pairing of spin manifolds and spin manifolds with singularities (spin bordism cycles and spin \(\mathbb{Z}_{k}\)-cycles), through the Conner–Floyd type isomorphism on \(KO\)-homology due to Hovey and Hopkins.
These \(\mathbb{Q}\)- and \(\mathbb{Q}/\mathbb{Z}\)-periods are proved to form a complete set of numerical invariants for Real \(K\)-theory in [Hu23].
The goal of this section is to show our Differential KO character will induce the complete set of numerical invariants on Real K theory described above.

\paragraph{Definition 5.1}
Let $h^{*}$ be a generalized cohomology theory of finite type (each $h^{i}(\mathrm{pt})$ finitely
generated). The exact functors $\operatorname{Hom}(-,\mathbb{Q})$ and
$\operatorname{Hom}(-,\mathbb{Q}/\mathbb{Z})$ send $h^{*}(-)$ to generalized homology theories,
which (by Brown representability) are represented by spectra we denote
$D_{\mathbb{Q}}h$ and $D_{\mathbb{Q}/\mathbb{Z}}h$. The quotient map
$\mathbb{Q}\twoheadrightarrow \mathbb{Q}/\mathbb{Z}$ induces a map of spectra
\[
D_{\mathbb{Q}}h \longrightarrow D_{\mathbb{Q}/\mathbb{Z}}h .
\]
the \textbf{Anderson dual spectrum} $Dh$ of $h$ is the homotopy fiber of this map.
The homology and cohomology theories represented by $Dh$ are called the
\emph{Anderson dual homology} $Dh_{*}$.

\paragraph{Theorem 5.2}[Hu 23]
Let $h^{*}$ be a generalized cohomology theory of finite type and let $D h_{*}$
be its Anderson dual. If $X$ is a finite CW–complex, then for every $i$ there is a
natural isomorphism
\[
h^{\,i}(X)\;\cong\;
\left[
\vcenter{\hbox{%
\begin{tikzcd}[column sep=4em,row sep=3em,baseline=(current bounding box.center)]
(D_{\mathbb{Q}} h)_{i}(X) \arrow[r] \arrow[d] &
\mathbb{Q} \arrow[d, "\mathrm{mod}\,\mathbb{Z}"] \\
(D_{\mathbb{Q}/\mathbb{Z}} h)_{i}(X) \arrow[r] &
\mathbb{Q}/\mathbb{Z}
\end{tikzcd}%
}}
\right],
\]

\begin{proof}

by [And69]
\begin{enumerate}
\item \(Dh_{i}(X;\mathbb{Q}) \cong \operatorname{Hom}\bigl(h^{i}(X),\mathbb{Q}\bigr),\) \emph{and}
\item \(Dh_{i}(X;\mathbb{Q}/\mathbb{Z}) \cong \operatorname{Hom}\bigl(h^{i}(X),\mathbb{Q}/\mathbb{Z}\bigr).\)
\end{enumerate}

Consider \(A\) be a finitely generated abelian group. Then the evaluation map
\[
A\longrightarrow
\left\{
\begin{tikzcd}[ampersand replacement=\&, column sep=large]
\Hom(A,\mathbb{Q}) \arrow[r] \arrow[d] \& \mathbb{Q} \arrow[d] \\
\Hom(A,\mathbb{Q}/\mathbb{Z}) \arrow[r] \& \mathbb{Q}/\mathbb{Z}
\end{tikzcd}
\right\}
\]

It suffices to show \(A=\mathbb{Z}_{n}\) and \(A=\mathbb{Z}\).
For \(A=\mathbb{Z}_{n}\) one has \(\Hom(\mathbb{Z}_{n},\mathbb{Q})=0\) and
\(\Hom(\mathbb{Z}_{n},\mathbb{Q}/\mathbb{Z})\cong \mathbb{Z}_{n}\), so the claim is immediate.
For \(A=\mathbb{Z}\), the assertion amounts to: every endomorphism
\(\varphi:\mathbb{Q}/\mathbb{Z}\to\mathbb{Q}/\mathbb{Z}\) that lifts to
\(\tilde\varphi:\mathbb{Q}\to\mathbb{Q}\) must be multiplication by some rational \(q\).
For \(\tilde\varphi\) to descend to \(\varphi\), one needs \(q\mathbb{Z}\subset\mathbb{Z}\),
hence \(q\in\mathbb{Z}\). Thus \(\varphi\) is multiplication by an integer, which proves the
isomorphism. 
\end{proof}

It was known that the Anderson Dual of $KO^i$ is symplectic K homology, which is $KO$ shifted by 4. we have the following corrollary. 
\paragraph{Corollary 5.3}
Let $X$ be a finite CW–complex. For every $i$ there is a natural isomorphism
\[
KO^{0}(X)\;\cong\;
\left[
\begin{tikzcd}[column sep=4.2em,row sep=3em,baseline=(current bounding box.center)]
KO_{-4}(X;\mathbb{Q}) \arrow[r] \arrow[d] &
\mathbb{Q} \arrow[d, "\mathrm{mod}\,\mathbb{Z}"] \\
KO_{-4}(X;\mathbb{Q}/\mathbb{Z}) \arrow[r] &
\mathbb{Q}/\mathbb{Z}
\end{tikzcd}
\right],
\]

With Corrollary 5.3, we immediately yield the following theorem
\paragraph{Theorem 5.4}
A cohomology class $[\omega]\in H^{4\bullet}(X,\mathbb{Q})$ is the Pontryagin character of a real vector bundle over $X$ if and only if for every closed spin even-dimensional manifold $M^{4k}$ mapping to $X$, $M \xrightarrow{\,f\,} X$, the quantity
\[
\int_{M} \widehat{A}(TM)\wedge f^{*}\omega
\]
is an integer; moreover, when $\dim M\equiv 4 \pmod{8}$, this value is even.
\begin{proof}
From the short exact sequence of coefficients
\[
0\longrightarrow \mathbb Z \xrightarrow{\;\iota\;} \mathbb Q
\xrightarrow{\;\pi\;} \mathbb Q/\mathbb Z \longrightarrow 0
\]
we obtain the long exact (Bockstein) sequence in real $K$–homology
\[
\cdots \to \KO_{4}(X)\xrightarrow{\ \iota_* \ } \KO_{4}(X;\mathbb Q)
\xrightarrow{\ \pi_* \ } \KO_{4}(X;\mathbb Q/\mathbb Z)
\xrightarrow{\ \beta\ } \KO_{3}(X)\to \cdots .
\]

\[
\begin{tikzcd}[column sep=large, row sep=large]
{} \arrow[r, "\beta"] &
KO_{4}(X)
  \arrow[r, "\otimes \mathbb{Q}"]
  \arrow[d, "C_{\mathbb Z}"'] &
KO_{4}(X)\otimes \mathbb{Q}
  \arrow[r]
  \arrow[d, "C_{\mathbb Q}"'] &
KO_{4}(X;\mathbb{Q}/\mathbb{Z})
  \arrow[r, "\beta"]
  \arrow[d, "C"'] &
KO_{3}(X;\mathbb{Z})
\\
0 \arrow[r] &
\mathbb{Z} \arrow[r] &
\mathbb{Q} \arrow[r] &
\mathbb{Q}/\mathbb{Z} \arrow[r] &
0
\end{tikzcd}
\]

Let $\omega\in\Omega^{4\bullet}(X)$ be a closed form satisfying the
$\Ahat$–integrality condition:
for every closed spin $4k$–manifold $M$ and smooth $f:M\to X$,
\[
\int_M \Ahat(TM)\wedge f^*\omega\ \in\ \mathbb Z,
\quad\text{and for }\dim M\equiv 4\ (\mathrm{mod}\ 8)\text{ the value is even.}
\]
By Conner and Floyd type isomorphism developed by Hopkins and Hovey[18]. 
\[
KO_{*}(X)\ \cong\
M\!Spin_{*}(X)\Big/\Big\langle\ [Q]\cdot x - \alpha(Q)\,x\ :\ Q\in M\!Spin_{*},\ x\in M\!Spin_{*}(X)\ \Big\rangle .
\]
Here \([Q]\cdot x\) is the \(M\!Spin_{*}\)-module action (external product \(Q\times -\) followed by projection to \(X\)), and \(\alpha(Q)\in KO_{*}\) is the \(KO\)-valued Dirac index/A hat genus of \(Q\). 
This defines a homomorphism
\[
C_{\mathbb Z}:\ \KO_{4}(X)\longrightarrow \mathbb Z,\qquad
C_{\mathbb Z}([f:M\to X])\ :=\ \int_M \Ahat(TM)\wedge f^*\omega .
\]
Tensoring with $\mathbb Q$ kills torsion in $\KO_{4}(X)$, so we obtain a
well-defined rational functional
\[
C_{\mathbb Q}:\ \KO_{4}(X;\mathbb Q)\cong \KO_{4}(X)\otimes\mathbb Q
\longrightarrow \mathbb Q ,
\]
supported only on the  spin $4k$–cycles.

Since $\mathbb Q/\mathbb Z$
is a divisible (hence injective) $\mathbb Z$–module, any homomorphism
defined on a subgroup extends across an inclusion (Baer criterion).
Applied to the short exact sequence above, the map $C_{\mathbb Q}$ 
 extends uniquely modulo $\mathbb Z$ to a homomorphism
\[
C:\ \KO_{4}(X;\mathbb Q/\mathbb Z)\longrightarrow \mathbb Q/\mathbb Z
\]
such that $C\circ \pi_*=\pi\circ C_{\mathbb Q}$, i.e. $C$ is a lift of
$C_{\mathbb Q}$ through the quotient $\mathbb Q\to \mathbb Q/\mathbb Z$.

By theorem 5.3 ,
a real $K$–theory cocycle (equivalently, a real vector bundle) is determined by such "liftable" $C$ . 

The other way follows from Atiyah-Singer index theory.
\end{proof}

\paragraph{Definition 5.5}
We recall the definition of $\wedge_{O}$ are forms that are cohomologous to the odd pontryagin character + exact.  
\
The odd Pontryagin character is the natural map
\(ph_{\mathrm{odd}}:KO^{-1}(X)\to \bigoplus_{k\ge0}H^{4k-1}(X;\mathbb{Q})\)
that commutes with suspension in \(KO\) and the Pontryagin character, i.e., the diagram commutes.
\[
\begin{tikzcd}
\widetilde{KO}^{0}(\Sigma X) \ar[r,"{\mathrm{ph}}"] \ar[d,"\Sigma^{-1}"',"\cong" sloped] &
\displaystyle \bigoplus_{k\ge0}\,\widetilde{H}^{4k}(\Sigma X;\mathbb{Q})
\ar[d,"\sigma","\cong" sloped] \\
\widetilde{KO}^{-1}(X) \ar[r,"{\mathrm{ph}_{\mathrm{odd}}}"'] &
\displaystyle \bigoplus_{k\ge0}\,\widetilde{H}^{4k-1}(X;\mathbb{Q})
\end{tikzcd}
\]

\noindent\textbf{Corrollary 5.6}
let $w$ be a closed form of $4k-1$ degrees on X. $w \in \wedge_{O}$ if the following holds:
for every closed Spin manifold $M$ of dimension $8K+3$ and $8K+7$,
and every smooth map $f:M\to X$, we have
\[
\int_{M}\widehat A(TM)\wedge f^{*}\omega \in
\begin{cases}
\mathbb Z, & \dim M\equiv 7 \pmod 8,\\[2pt]
2\mathbb Z, & \dim M\equiv 3 \pmod 8.
\end{cases}
\]

\begin{proof}
The corollary follows from applying the proof of theorem 5.3.1 with suspension isomorphism to \(KO^0\), and \(KO_{4}(X;\mathbb{Q}/\mathbb{Z})\).

\end{proof}

\paragraph{Definition 5.7}
A \emph{\(\mathbb{Z}_{k}\)-manifold} is a pair \((M,\beta M)\) consisting of
\begin{itemize}
  \item an oriented smooth manifold \(M\) with boundary, and
  \item a closed oriented manifold \(\beta M\),
\end{itemize}
such that the boundary decomposes as a disjoint union of \(k\) (labelled) copies of
\(\beta M\):
\[
\partial M \;\cong\; \bigsqcup_{j=1}^{k} (\beta M)_{j}.
\]
Let \(\overline M\) be the quotient obtained by gluing each boundary component of \(M\)
to \(\beta M\) by the given labels; we also call \(\overline M\) a \(\mathbb{Z}_{k}\)-manifold.
The manifold \(\beta M\) is called the \emph{Bockstein} of \(\overline M\).

If \((N,\beta N)\) is another \(\mathbb{Z}_{k}\)-manifold, we say that
\((M,\beta M)\) is \emph{\(\mathbb{Z}_{k}\)-cobordant} to \((N,\beta N)\) if there exists
a cobordism \(W\) from \(\beta M\) to \(\beta N\) (so
\(\partial W=\beta W \,\sqcup\, (-\beta N)\)) such that the glued manifold
\[
M \;\cup_{\partial M}\; \bigsqcup_{j=1}^{k} W \;\cup_{-\partial N}\; (-N)
\]
is a boundary.  A \(\mathbb{Z}_{k}\)-manifold is a \emph{\(\mathbb{Z}_{k}\)-boundary}
if it is \(\mathbb{Z}_{k}\)-cobordant to the empty \(\mathbb{Z}_{k}\)-manifold.

\begin{figure}[htbp]
  \centering
  \includegraphics[width=0.6\textwidth]{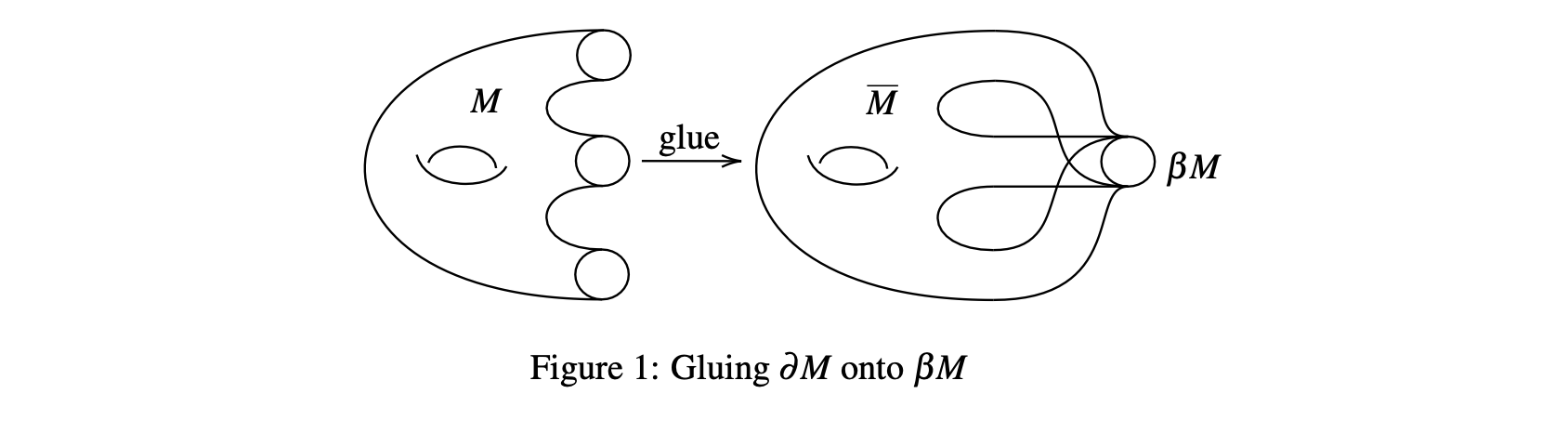}
  \caption{Example of a \(\mathbb{Z}_{3}\)-manifold.}
  \label{fig:zk}
\end{figure}

We write
\[
\Omega_*^{\mathrm{Spin}}(X;\mathbb{Z}_k)
\]
to denote Spin bordism cycles over $X$ equipped with a $\mathbb{Z}_k$–manifold structure.

\paragraph{Definition 5.8}
Let $E$ be a real vector bundle.

We define the angle invariant above on a \(\mathbb{Z}_{k}\)-spin manifold is exactly \(\mathrm{inv}_{\mathbb{Z}_{k}}\), where
\[
\mathrm{inv}_{\mathbb{Z}_{k}}(E):
\Omega^{\mathrm{Spin}}_{*}\!\left(X;\mathbb{Z}_{k}\right)
\xrightarrow{\ \widehat{A}\ }
KO_{*}\!\left(X;\mathbb{Z}_{k}\right)
\xrightarrow{\ -\!\setminus E\ }
KO_{*}\!\left(\mathrm{pt};\mathbb{Z}_{k}\right).
\]

Passing to the direct limit over \(k\),
\[
\varinjlim_{k}\,\mathrm{inv}_{\mathbb{Z}_{k}}(E)
\;=\;
\mathrm{inv}_{\mathbb{Q}/\mathbb{Z}}(E),
\qquad
\text{since } \varinjlim_{k}\mathbb{Z}_{k}\cong \mathbb{Q}/\mathbb{Z}.
\]

In fact, one can apply this pairing to $\mathbb{Q}$.
\[
\mathrm{inv}_{\mathbb{Q}}(E):
\Omega^{\mathrm{Spin}}_{*}\!\left(X;\mathbb{Q}\right)
\xrightarrow{\ \widehat{A}\ }
KO_{*}\!\left(X;\mathbb{Q}\right)
\xrightarrow{\ -\!\setminus E\ }
KO_{*}\!\left(\mathrm{pt};\mathbb{Q}\right).
\]

\paragraph{Theorem 5.9}
Differential \(KO\)-character induces complete set of numerical invariants determining a real vector bundle up to stably equivalence.

\begin{proof}

If $(\overline{M},f)$ is a \emph{spin} $\mathbb{Z}_k$-cycle in $X$ of dimension $n \equiv 2,3 \pmod{8}$, then the evaluation of $\mathrm{inv}_{\mathbb{Q}/\mathbb{Z}}(E)$ on $(\overline{M},f)$ , via the Bockstein isomorphism \[
\mathrm{KO}_n(\mathrm{pt};\mathbb{Q}/\mathbb{Z}) \cong \mathrm{KO}_{n-1}(\mathrm{pt}),
\] becomes the $\mathbb{Z}_2$ index of spin dirac operator on the bockstein twisted by the pull back bundle, which is by definition $\hat{\kappa}_{(E,\nabla)}(\beta M)$.The Bockstein $\beta M$ has dimension $8K+1$ or $8K+2$.

\[
\begin{tikzcd}[column sep=large,row sep=large]
\Omega^{\mathrm{Spin}}_{n}\!\left(X;\,\mathbb{Q}/\mathbb{Z}\right)
  \arrow[r, "\overline{M}\mapsto \beta M", "\mathrm{Bockstein}"']
  \arrow[d, "\mathrm{inv}_{\mathbb{Q}/\mathbb{Z}}(E)"']
& \Omega^{\mathrm{Spin}}_{\,n-1}(X)
  \arrow[d, "\widehat{A}\text{ slant with }E\;=\;\mathrm{ind}\!\big(D_{\mathrm{spin}}\otimes f^{*}E\big)"] \\[2pt]
KO_{n}\!\left(\mathrm{pt};\,\mathbb{Q}/\mathbb{Z}\right)
  \arrow[r, "\mathrm{Bockstein}"', "\;\cong\;"{sloped}]
& KO_{\,n-1}(\mathrm{pt})
\end{tikzcd}
\qquad
\text{with } n\in\{8k+2,\,8k+3\}.
\]

We note the dimensions \(8K\), and \(8K+4\).
For a \(\mathbb{Z}_{k}\)-manifold represented by the pair \((\overline{M},\,\beta M)\) and a map \(f:\overline{M}\to X\), 
\[
\mathrm{inv}_{\mathbb{Z}_{k}}(\overline{M},\beta M)
\;:=\;
\Biggl[
\int_{\overline{M}}\widehat{A}(T\overline{M})\wedge \operatorname{ph}\!\bigl(f^{*}E,f^{*}\nabla\bigr)
\;-\;
k\,\hat{\kappa}_{(E,\nabla)}(\beta M)
\Biggr]\ \bmod k.
\]
Here \(\hat{\kappa}_{(E,\nabla)}(\beta M)\in \mathbb{R}/\mathbb{Z}\) is the value of the \(KO\)-character on the Bockstein boundary \(\beta M\).
This is true because by Mod K index theorem in [FM88] on $4k-1$ spin manifold $\mathrm{inv}_{\mathbb{Z}_{k}}$ admits an analytic index.  

The rational information of the bundle is recovered from the Pontryagin  form associated to the KO-character, ie \[
\operatorname{Hom}\!\Bigl(\Omega^{\mathrm{Spin}}_{4*}(X)\big/\mathrm{tors},\,\mathbb{Z}\Bigr),
\]
given by
\[
[f\colon M^{4k}\!\to X]\ \longmapsto\ \int_{M}\widehat{A}(TM)\wedge f^{*}\!\operatorname{ph}(E),
\]

We thus conclude $\hat{\kappa}_{(E,\nabla)}(\beta M)$ induces $\mathrm{inv}_{\mathbb{Q}/\mathbb{Z}}(E)$ and $\mathrm{inv}_{\mathbb{Q}}(E)  on (\overline{M},f)$ 

To finish off the proof, we now procede with Hu's argument.

The maps \(\mathrm{inv}^{}_{\mathbb{Q}}\) and \(\mathrm{inv}^{}_{\mathbb{Q}/\mathbb{Z}}\) are compatible in the sense that,
for any real vector bundle \(E\) over \(X\), the diagram
\[
\begin{tikzcd}[column sep=large]
\Omega^{\mathrm{Spin}}_{*}(X;\mathbb{Q})
  \arrow[r,"\mathrm{inv}^{}_{\mathbb{Q}}(E)"] \arrow[d, "\mathrm{mod}\ \mathbb{Z}"']
& KO_{*}(\mathrm{pt};\mathbb{Q}) \arrow[d, "\mathrm{mod}\ \mathbb{Z}"] \\
\Omega^{\mathrm{Spin}}_{*}(X;\mathbb{Q}/\mathbb{Z})
  \arrow[r,"\mathrm{inv}^{}_{\mathbb{Q}/\mathbb{Z}}(E)"']
& KO_{*}(\mathrm{pt};\mathbb{Q}/\mathbb{Z})
\tag{1}
\end{tikzcd}
\]
commutes, because slant products are compatible with change of coefficients.
Consider the map
\[
\mathrm{inv}^{}:\ KO(X)\longrightarrow
\left\{
\begin{tikzcd}[column sep=large]
\Omega^{\mathrm{Spin}}_{*}(X;\mathbb{Q})
  \arrow[r] \arrow[d, "\mathrm{mod}\ \mathbb{Z}"']
& KO_{*}(\mathrm{pt};\mathbb{Q}) \arrow[d, "\mathrm{mod}\ \mathbb{Z}"] \\
\Omega^{\mathrm{Spin}}_{*}(X;\mathbb{Q}/\mathbb{Z})
  \arrow[r]
& KO_{*}(\mathrm{pt};\mathbb{Q}/\mathbb{Z})
\end{tikzcd}
\right\},
\]
which sends \(E\) to the diagram \((1)\). Here the target is the group of commutative
diagrams of this form.

Both \(\mathrm{inv}^{}_{\mathbb{Q}}(E)\) and \(\mathrm{inv}^{}_{\mathbb{Q}/\mathbb{Z}}(E)\) are equivariant with respect to
the \(\widehat A\)-orientation \(\widehat A:\Omega^{\mathrm{Spin}}_{*}(\mathrm{pt})\to KO_{*}(\mathrm{pt})\).
Therefore \(\mathrm{inv}^{}\) lands in the subgroup of diagrams that are equivariant for
\(\widehat A\). By tensor–hom adjunction, this subgroup is isomorphic to the group
\[
\left\{
\begin{tikzcd}[column sep=large]
\Omega^{\mathrm{Spin}}_{*}(X;\mathbb{Q}) \arrow[r] \arrow[d, "\mathrm{mod}\ \mathbb{Z}"']
& KO_{*}(\mathrm{pt};\mathbb{Q}) \arrow[d, "\mathrm{mod}\ \mathbb{Z}"] \\
\Omega^{\mathrm{Spin}}_{*}(X;\mathbb{Q}/\mathbb{Z}) \arrow[r]
& KO_{*}(\mathrm{pt};\mathbb{Q}/\mathbb{Z})
\end{tikzcd}
\right\}_{KO_{*}(\mathrm{pt})}
\tag{2}
\]
of commutative diagrams of \(KO_{*}(\mathrm{pt})\)-modules of the above form, where
\(\Omega^{\mathrm{Spin}}_{*}(X;\Lambda)\) denotes
\(\Omega^{\mathrm{Spin}}_{*}(X)\otimes_{\Omega^{\mathrm{Spin}}_{*}(\mathrm{pt})}KO_{*}(\mathrm{pt})\) for
\(\Lambda=\mathbb{Q},\mathbb{Q}/\mathbb{Z}\).
By hopkins and hovey's conner and floyd type isomorphism for real K theory, the natural transformation \(\widehat A\) induces split surjection
\(\Omega^{\mathrm{Spin}}_{*}(X;\Lambda)\twoheadrightarrow KO_{*}(X;\Lambda)\).
Pulling back along these split surjections yields an embedding
\(i\) of \((2)\), as a direct summand, into
\[
\left\{
\begin{tikzcd}[column sep=large]
KO_{*}(X;\mathbb{Q}) \arrow[r] \arrow[d, "\mathrm{mod}\ \mathbb{Z}"']
& KO_{*}(\mathrm{pt};\mathbb{Q}) \arrow[d, "\mathrm{mod}\ \mathbb{Z}"] \\
KO_{*}(X;\mathbb{Q}/\mathbb{Z}) \arrow[r]
& KO_{*}(\mathrm{pt};\mathbb{Q}/\mathbb{Z})
\end{tikzcd}
\right\}_{KO_{*}(\mathrm{pt})}.
\tag{3}
\]

now projecting (3) to degree -4, and apply corrolary 5.3. 

\end{proof}

\section{Isomorphism between Differential KO theory and Differential KO character}

\vspace{.5cm}
\begin{center}
\setlength{\unitlength}{0.5cm}
\begin{picture}(24,16)(2.5,0)
\thicklines
\put(5,1){$0$}
\put(20.5,1){$0$}

\put(6,2){\vector(1,1){1.5}}
\put(18,3.5){\vector(1,-1){1.5}}

\put(8,4.5){$\wedge^{4k-1}/\wedge_{O}$}
\put(12,4.5){\vector(1,0){2.5}}
\put(16.5,4.5){$\wedge_{BO}$}
\put(13,5){\small{$d$}}
 
\put(6.5,7.5){\vector(1,-1){1.5}}
\put(10.5,7){\small{$i_{2}$}}
\put(10.5,6){\vector(1,1){1.5}}
\put(14.5,7.5){\vector(1,-1){1.5}}
\put(15.5,7){\small{$\widehat{ph}$}}
\put(18.5,6){\vector(1,1){1.5}}

\put(0,8){1)}
\put(3,8){$H^{4k-1}(\mathbb{R})$}
\put(12.75,8){$\hat{KO}$}
\put(20,8){$H^{4k}(\mathbb{R})$}

\put(6,9.5){\vector(1,1){1.5}}
\put(10.5,11){\vector(1,-1){1.5}}
\put(11.5,10.5){\small{$i_{1}$}}
\put(14,10.5){\small{$\delta_{2}$}}
\put(14.0,9.5){\vector(1,1){1.5}}
\put(18,11){\vector(1,-1){1.5}}

\put(7,12){$KO^{-1}(\mathbb{R}/\mathbb{Z})$}
\put(12,12){\vector(1,0){2.5}}
\put(15.5,12){$KO^{0}(\mathbb{Z})$}

\put(5.5,14.5){\vector(1,-1){1.5}}
\put(18,13){\vector(1,1){1.5}}

\put(4.5,15){$0$}
\put(20,15){$0$}
\end{picture}
\end{center}

\vspace{.5cm}
\begin{center}
\setlength{\unitlength}{0.5cm}
\begin{picture}(24,16)(2.5,0)
\thicklines
\put(5,1){$0$}
\put(20.5,1){$0$}

\put(6,2){\vector(1,1){1.5}}
\put(18,3.5){\vector(1,-1){1.5}}

\put(4.5,4.5){$\wedge^{4k-1}/\wedge_{4k-1-integrality}$}
\put(12,4.5){\vector(1,0){2.5}}
\put(16.5,4.5){$\wedge_{integrality}$}
\put(13,5){\small{$d$}}
 
\put(6.5,7.5){\vector(1,-1){1.5}}
\put(10.5,7){\small{$i_{2}$}}
\put(10.5,6){\vector(1,1){1.5}}
\put(14.5,7.5){\vector(1,-1){1.5}}
\put(15.5,7){\small{$\delta_{1}$}}
\put(18.5,6){\vector(1,1){1.5}}

\put(0,8){2)}  
\put(3,8){$
\Hom\!\big(KO_{3}(X),\,\mathbb{R}\big)
$}
\put(12,8){$\hat{KOCH}$}
\put(20,8){$
\Hom\!\big(KO_{4}(X),\,\mathbb{R}\big)
$}

\put(6,9.5){\vector(1,1){1.5}}
\put(10.5,11){\vector(1,-1){1.5}}
\put(11.5,10.5){\small{$i_{1}$}}
\put(14,10.5){\small{$\delta_{2}$}}
\put(14.0,9.5){\vector(1,1){1.5}}
\put(18,11){\vector(1,-1){1.5}}

\put(3,12){$
\Hom\!\big(KO_{3}(X),\,\mathbb{R}/\mathbb{Z}\big)
$}
\put(12,12){\vector(1,0){2.5}}
\put(15.5,12){$\mathrm{RationalHom}\!\big(KO_{4}(X;\mathbb{Q}/\mathbb{Z}),\,\mathbb{Q}/\mathbb{Z}\big)$}

\put(5.5,14.5){\vector(1,-1){1.5}}
\put(18,13){\vector(1,1){1.5}}

\put(4.5,15){$0$}
\put(20,15){$0$}
\end{picture}
\end{center}

$ KO_{*}(X)$ is interpreted as the Conner and floyd isomorphism for KO homology developed in [18] by Hovey and Hopkins.

$\Omega^{\mathrm{Spin}}_{*}(X)\otimes_{\Omega^{\mathrm{Spin}}_{*}} KO_{*}\ \xrightarrow{\ \cong\ }\ KO_{*}(X).$
\medskip

\(\wedge_{\mathrm{integrality}}(X)\) and
\(\wedge_{4k-1\text{-integrality}}(X)\) denote the closed forms on \(X\)
satisfying the following integrality conditions.

\smallskip
\emph{(Even degrees).} A closed 4k total form \(C\) belongs to
\(\wedge_{\mathrm{integrality}}(X)\) iff for every closed Spin manifold
\(M^{4k}\) and every smooth map \(f:M\to X\),
\[
\int_{M} f^{*}C \wedge \widehat A\!\big(\nabla^{TM}\big)\ \in\ \mathbb{Z},
\]
and moreover this integer is \emph{even} when \(k\) is even.

\smallskip
\emph{(Odd degrees \(4k-1\)).} A closed total \((4k-1)\)-form \(\omega\) belongs to
\(\wedge_{4k-1\text{-integrality}}(X)\) iff for every closed Spin manifold
\(M^{4k-1}\) and every smooth map \(g:M\to X\),
\[
\int_{M} g^{*}\omega \wedge \widehat A\!\big(\nabla^{TS}\big)\ \in\ \mathbb{Z}.
\]
and  this integer is \emph{even} when \(k\) is even.

\noindent\emph{Note.} By Theorem~5.4 and Corollary~5.6 we have the identifications
\[
\wedge_{\mathrm{integrality}}(X)\;=\;\wedge_{BO}(X),
\qquad
\wedge_{4k-1\text{-integrality}}(X)\;=\;\wedge_{\circ}(X).
\]

\noindent\textbf{Theorem 6.1}
There is a natural isomorphism with resepct to the hexagon
\[
\widehat{KO}(X) \xrightarrow{\ \cong\ } \widehat{\mathrm{KOCH}}(X),
\]
sending a real vector bundle with metric connection \((E,\nabla)\) to its associated differential \(KO\)-character.

The entire section aims to prove theorem 6.1. 

\paragraph{Uniqueness philosophy.}
The hexagon serves as a \emph{criterion for uniqueness}: once the hexagon
diagram (with exact rows and natural structure maps) is fixed, any theory
fitting into the middle and making the diagram commute is canonically
isomorphic to any other such theory. We apply this principle to identify
\emph{KO-characters} with \emph{differential \(KO\)-theory}, yielding
a natural isomorphism \( \widehat{KOCh}(X) \cong \widehat{KO}^{0}(X) \).
For background and comparison, see [12],[32],[33]\par

We now define the 4 diagonal maps appearing in diagram (2).

\medskip

\noindent\textbf{Definition 6.1}

Define
\[
i_{2}:\ \frac{\wedge^{\mathrm{4k-1}}(X)}{\wedge^{\mathrm{4k-1}}_{\mathrm{integrality}}(X)}
\longrightarrow \widehat{\mathrm{KOCH}}(X)
\]
as follows. For \(\phi\in\Lambda^{4k-1}\), by isomorphism in theorem 2.9, we obtain a stably trivial bundle with
connection \((E,\nabla^{E})=\widehat{CS}^{-1}(\phi)\). Then \(i_{2}(\phi)\) is the
differential \(KO\)-character whose value on an enriched Spin cycle \(f:M\to X\) is
\[
i_{2}(\phi)(f)\ :=\
\bar\eta\!\left(D_{M}\otimes f^{*}(E,\nabla^{E})\right)
\;-\;\mathrm{rk}(E)\,\bar\eta(D_{M})
\quad (\mathrm{mod}\ \mathbb{Z}),
\]
i.e. the rho invariant.

\noindent\textbf{Proposition 6.2}
The assignment \(i_{2}\) is well defined.
Equivalently, if \[
\phi \in \Lambda^{4k-1}_{\mathrm{integrality}}(X)
\] , ie
\(\widehat{CS}^{-1}(\phi)\) is  stably trivial and stably flat then for every enriched Spin cycle
\(f:M\to X\),
\[
i_{2}(\phi)(f)\;=\;
0\quad(\mathrm{mod}\ \mathbb{Z}).
\]

\begin{proof}
Let $\{\nabla_t\}_{t\in[0,1]}$ be a path on $E\oplus\underline{\mathbb R}^{\,k}$ with
$\nabla_0=\nabla\oplus d$ and $\nabla_1=\nabla_{\mathrm{flat}}$ flat. By additivity of reduced eta,
\[
i_2(\phi)(f)\ \equiv\ 
\bar\eta\!\big(D_M^{\,f^*(\nabla\oplus d)}\big)-\bar\eta\!\big(D_M^{\,f^*\nabla_{\mathrm{flat}}}\big)
\pmod{\mathbb Z}.
\]
The APS variation formula along $\{\nabla_t\}$ gives
\[
\bar\eta\!\big(D_M^{\,f^*(\nabla\oplus d)}\big)-\bar\eta\!\big(D_M^{\,f^*\nabla_{\mathrm{fl}}}\big)
\ \equiv\ \int_M \widehat A(\nabla^{TM})\wedge f^*\operatorname{CS}\!\big(\nabla_{\mathrm{flat}},\,\nabla\oplus d\big)
\pmod{\mathbb Z}.
\]
Clutch $\{\nabla_t\}$ to a bundle $\mathcal E\to M\times S^1$. Transgression yields
\[
\int_M \widehat A(\nabla^{TM})\wedge f^*\operatorname{CS}\!\big(\nabla_{\mathrm{fl}},\,\nabla\oplus d\big)
=\int_{M\times S^1}\widehat A\!\big(\nabla^{T(M\times S^1)}\big)\wedge \mathrm{ph}(\mathcal E).
\]
By the Atiyah–Singer index theorem, the right-hand side equals
\(\mathrm{index}\,D_{M\times S^1}^{\,\mathcal E}\in\mathbb Z\).
Therefore
\[
i_2(\phi)(f)\ \equiv\ \mathrm{index}\,D_{M\times S^1}^{\,\mathcal E}\ \equiv\ 0\quad \text{in }\mathbb R/\mathbb Z,
\]
as claimed.
\end{proof}

\noindent\textbf{Corollary 6.3}

\[
i_{2}(\phi)\big(f\big)
\;=\;
\int_{M}\widehat{A}(TM,\nabla^{TM})\wedge f^{*}\phi
\ \ (\mathrm{mod}\ \mathbb{Z}),
\]

\noindent\textbf{Definition 6.4}

\[
\mathrm{RationalHom}\big(KO_{4}(X;\mathbb{Q}/\mathbb{Z}),\mathbb{Q}/\mathbb{Z}\big)
\]

we mean the subset of \(Hom\big(KO_{4}(X;\mathbb{Q}/\mathbb{Z}),\mathbb{Q}/\mathbb{Z}\big)\)
consisting of those functionals that \emph{lift} through the coefficient
quotient \(\pi:\mathbb{Q}\twoheadrightarrow\mathbb{Q}/\mathbb{Z}\); namely
\[
\mathrm{RationalHom}\big(KO_{4}(X;\mathbb{Q}/\mathbb{Z}),\mathbb{Q}/\mathbb{Z}\big)
=\Big\{\phi\ \Big|\ \exists\,\widetilde\phi:\ KO_{4}(X)\otimes\mathbb{Q}\to\mathbb{Q}
\ \text{ with }\ \pi\circ\widetilde\phi=\phi\circ\pi_{*}\Big\}.
\]

with such a lift \(\widetilde\phi\) it represents a \emph{real \(K\)-theory cocycle}.

\textbf{Definition 6.5}
$\delta_{2}$ denotes the map from differential $KO$-characters to
\[
\mathrm{RationalHom}\big(KO_{4}(X;\mathbb{Q}/\mathbb{Z}),\mathbb{Q}/\mathbb{Z}\big).
\]
The map is described in Theorem~5.9.

\textbf{Proposition 6.6 }
$\delta_{2}$ is surjective.

\begin{proof}
    This is equivalent to theorem 5.9
\end{proof}

\textbf  {Definition 6.7}
The map $\delta_{1}$ extracts the \emph{associated integrality form} of a
differential KO–character when one varies by cobordism.

\textbf{Proposition 6.8  }
$\delta_{1}$ is surjective. 
\begin{proof}
From theorem 5.4 we see $\wedge_{BO}=\wedge_{integrality}$ \par
let $\mu \in \wedge_{BO}(X)$. By definition, there exist $\mathcal{V}\in \mathrm{Struct}(X)$ and $\theta \in \wedge^{4k-1}$ such that
\[
\mu = \operatorname{ph}(\mathcal{V}) + d\theta.
\]
Thus,
\[
\mu = \operatorname{ph}\!\bigl(\mathcal{V} + i_{2}(\{\theta\})\bigr).
\]

\end{proof}

\textbf{Definition 6.9}
Let \(i_{1}\) denote the inclusion of \emph{flat} differential \(KO\)-characters, i.e., those whose associated integrality form is \(0\).

\medskip
\textbf {Proposition 6.10} 
Two hexagon at entries 1,3,6,7,5,2 are cannonically isomorphic. 

\vspace{.5cm}
\begin{center}
\setlength{\unitlength}{0.5cm}
\begin{picture}(24,16)(2.5,0)
\thicklines
\put(5,1){$0$}
\put(20.5,1){$0$}

\put(6,2){\vector(1,1){1.5}}
\put(18,3.5){\vector(1,-1){1.5}}

\put(9,4.5){$\mathbf{1}$}
\put(12,4.5){\vector(1,0){2.5}}
\put(16.5,4.5){$\mathbf{2}$}
\put(13,5){\small{$d$}}
 
\put(6.5,7.5){\vector(1,-1){1.5}}
\put(10.5,7){\small{$i_{2}$}}
\put(10.5,6){\vector(1,1){1.5}}
\put(14.5,7.5){\vector(1,-1){1.5}}
\put(15.5,7){\small{$\delta_{1}$}}
\put(18.5,6){\vector(1,1){1.5}}

\put(0,8){3)}  
\put(6,8){$\mathbf{3}$}
\put(13,8){$\mathbf{4}$}
\put(20,8){$\mathbf{5}$}

\put(6,9.5){\vector(1,1){1.5}}
\put(10.5,11){\vector(1,-1){1.5}}
\put(11.5,10.5){\small{$i_{1}$}}
\put(14,10.5){\small{$\delta_{2}$}}
\put(14.0,9.5){\vector(1,1){1.5}}
\put(18,11){\vector(1,-1){1.5}}

\put(9,12){$\mathbf{6}$}
\put(12,12){\vector(1,0){2.5}}
\put(16.5,12){$\mathbf{7}$}

\put(5.5,14.5){\vector(1,-1){1.5}}
\put(18,13){\vector(1,1){1.5}}

\put(4.5,15){$0$}
\put(20,15){$0$}
\end{picture}
\end{center}

\begin{proof}
 
\[
KO^{8k-1}(X;\mathbb{R}/\mathbb{Z}) \;\cong\;
\operatorname{Hom}\!\big(KO_{3}(X),\,\mathbb{R}/\mathbb{Z}\big),
\]
\[
H^{4k-1}(X;\mathbb{R}) \;\cong\;
\operatorname{Hom}\!\big(KO_{3}(X),\,\mathbb{R}\big),
\]
\[
H^{4k}(X;\mathbb{R}) \;\cong\;
\operatorname{Hom}\!\big(KO_{4}(X),\,\mathbb{R}\big),
\]

These three cannonical isomorphism follow from Universal Coefficient theorem on Anderson Duality of KO.

\[
KO^{8k}(X;\mathbb{Z}) \;\cong\;
\textit{Rational }\operatorname{Hom}\!\big(KO_{4}(X;\mathbb{Q}/\mathbb{Z}),\,\mathbb{Q}/\mathbb{Z}\big),
\]

Follows from Corollary 5.3

\[
\frac{\wedge^{4k-1}(X)}{\wedge^{4k-1}_{\mathrm{integrality}}(X)}
\;\cong\;
\frac{\wedge^{4k-1}(X)}{\wedge_{o}^{4k-1}(X)},
\qquad
\wedge_{BO}(X) \;\cong\; \wedge_{\mathrm{integrality}}(X).
\]
These two isomorphisms follow from theorem 5.4 and corollary 5.6.

\end{proof}

\noindent\textbf{Proposition 6.11 }

The following sequence is short exact
\[
0 \longrightarrow \operatorname{Hom}\!\big(KO_{3}(X),\mathbb{R}/\mathbb{Z}\big)
  \longrightarrow \widehat{\mathrm{KOCH}}(X)
  \xrightarrow{\ \delta_{1}\ } \wedge_{\mathrm{integrality}}
  \longrightarrow 0.
\]

\begin{proof}

By the Conner–Floyd isomorphism, we have a decomposition
\[
\operatorname{Hom}\!\bigl(KO_{3}(X),\mathbb{R}/\mathbb{Z}\bigr)
\;\cong\;
\left(
\operatorname{Hom}\!\bigl(\Omega^{\mathrm{Spin}}_{8K+1}(X),\mathbb{Z}_{2}\bigr)
\oplus
\operatorname{Hom}\!\bigl(\Omega^{\mathrm{Spin}}_{8K+2}(X),\mathbb{Z}_{2}\bigr)
\oplus
\operatorname{Hom}\!\bigl(\Omega^{\mathrm{Spin}}_{8K+3}(X),\mathbb{R}/\mathbb{Z}\bigr)
\right.
\]
\[
\left.
\oplus\,
\operatorname{Hom}\!\bigl(\Omega^{\mathrm{Spin}}_{8K+7}(X),\mathbb{R}/\mathbb{Z}\bigr)
\right)\Big/\!\!\sim
\]

where \(\sim\) imposes the equivalence of Conner and floyd isomorphism.

Hence, whenever a homomorphism is well defined on \(KO_{3}(X)\), it is a cobordism
invariant on spin bordism. Varying by a cobordism then yields value \(0\), which
implies the integral form of the associated differential KO character to be 0. Consequently,
\[
\ker(\delta_{1}) \;=\; \operatorname{im}(i_{1}).
\]

\end{proof}

\textbf{Proposition 6.12 }

\[
0\;\longrightarrow\;
\frac{\wedge^{\mathrm{4k-1}}(X)}{\wedge^{\mathrm{4k-1}}_{\mathrm{int}}(X)}
\xrightarrow{\ i_{2}\ }\widehat{\mathrm{KOCH}}(X)
\xrightarrow{\ \delta_{2}\ }
\mathrm{RationalHom}\!\big(KO_{4}(X;\mathbb{Q}/\mathbb{Z}),\mathbb{Q}/\mathbb{Z}\big)
\;\longrightarrow 0
\]

is short exact. 

\begin{proof}

\[
\text{We first show } 
\ker \delta_{2} \subseteq \operatorname{im} i_{2}.
\]

A zero element in real $K$-theory corresponds to a stably trivial bundle. 
In terms of 
\[
\mathrm{RationalHom}\!\big(KO_{4}(X;\mathbb{Q}/\mathbb{Z}),\mathbb{Q}/\mathbb{Z}\big),
\]
this means that all numerical invariants of the bundle vanish under pairing with every
$\mathbb{Z}_{k}$-Spin cycle in degrees $8K, 8K+2, 8K+3,$ and $8K+4$.

For a stably trivial real vector bundle \(E \to X\), the Pontryagin character
form is exact in positive degrees: for each \(k \ge 1\) there exists a form
\(\omega_{4k-1} \in \Omega^{4k-1}(X)\) such that
\[
d\omega_{4k-1} = \mathrm{ph}_{4k}(E,\nabla).
\]
In particular, the Pontryagin character form of \(E\) is exact in positive degrees.

For \(\omega\in\wedge^{\mathrm{4k-1}}(X)\) and a \(\mathbb{Z}_{k}\)-bordism class
\([f:M\to X]\) with Bockstein \(\beta M\) and
a bordism \(W\) satisfying \(\partial W=k\cdot \beta M\), we have by definition
\[
\delta_{2}\big(i_{2}(\omega)\big)\big([f]\big)
=\ \int_{W}\widehat A(W)\wedge f^{*}ph(E)
 \;-\; k\, i_{2}(\omega)(\beta M)\
\]
But
\[
k\, i_{2}(\omega)(\beta M)
=\int_{\partial W}\widehat A(T\partial W)\wedge f^{*}\omega
=\int_{W} d\!\big(\widehat A(TW)\wedge f^{*}\omega\big)
=\int_{W}\widehat A(TW)\wedge f^{*}(d\omega),
\]
we notice $dw=ph(E)$ by construction. 

Now we show 

\[\operatorname{im} i_{2} \subseteq \ker \delta_{2}.\]

This follows from the isomorphism
\  
\[
\widehat{CS}_{\!Ph}:\ \mathrm{Struct}_{\mathrm{ST}}(X)/\mathrm{Struct}_{\mathrm{SF}}(X)\ \cong\ \Lambda^{\mathrm{4k-1}}(X)/\Lambda_{O}(X) \cong\ \Lambda^{\mathrm{4k-1}}(X)/\Lambda_{integrality-4k-1}(X) .
\]

\end{proof}

The top row and bottom row of the hexagon comes from the bockestein sequence and de rham isomorphism respectively.  

\textbf{Remark.} We make a remark about the geometric meaning of the bockstein map. We note that the kernel of the natural map is the torsion part of Real K theory, i.e.
\[
\ker\!\left[\mathrm{RationalHom}\!\big(KO_{4}(X;\mathbb{Q}/\mathbb{Z}),\mathbb{Q}/\mathbb{Z}\big)\longrightarrow
\operatorname{Hom}\!\big(KO_{4}(X),\mathbb{R}\big)\right]
\;=\; \operatorname{tors}\, KO^{0}(X).
\]

 decomposes as
\[
\left(
\operatorname{Hom}\!\big(\operatorname{tors}\,\Omega^{\mathrm{Spin}}_{8k+2}(X),\mathbb{Z}_{2}\big)
\;\oplus\;
\operatorname{Hom}\!\big(\operatorname{tors}\,\Omega^{\mathrm{Spin}}_{8k+1}(X),\mathbb{Z}_{2}\big)
\;\oplus\;
\operatorname{Hom}\!\big(\operatorname{tors}\,\Omega^{\mathrm{Spin}}_{4k-1}(X),\mathbb{Q}/\mathbb{Z}\big)
\right)\Big/\!\!\sim
\]

Thus the bockstein map \(=\) restricting \(f:\operatorname{Hom}(KO_{3}(X),\mathbb{R}/\mathbb{Z})\) on \(\operatorname{tor} KO_{3}(X)\). The concrete index on \(\operatorname{tor} KO_{3}(X)\), represented by $\beta M$,  is obtained by lifting $\beta M$ (since it is represented by a torsion spin cycle) to a $\mathbb{Z}_k$ manifold and then evaluate with $\mathrm{inv}^{}_{\mathbb{Q}/\mathbb{Z}}(E)$ in definition 5.8. The resulting value is independent for the chocie of lifting, since one can verify $\mathrm{inv}^{}_{\mathbb{Q}/\mathbb{Z}}(E)$ is a $\mathbb{Z}_k$  cobordism invariant.

\medskip

Combining the propositions above and the uniqueness philosophy, we had thus proven theorem 6.1 .
\section{Family Index theorem in Differential KO theorem }

In this section we state the family index theorem for a Riemannian submersion
\(\pi:X\to B\) whose fibers are closed \(8k\)-dimensional \emph{Spin} manifolds,
in the setting of differential \(KO\)-theory for both the structured–bundle
model and the differential \(KO\)-character model. The fact that these two
pushforwards agree under the isomorphism described in Theorem 6.1 is precisely
the Bismut–Cheeger adiabatic limit theorem.

We begin by defining the pushforward in the Structured–bundle framework.
Via Proposition~3.2, this pushforward is obtained by translating the
Freed–Lott model to structured bundles, and hence all constructions and
identities below are the direct images of the Freed–Lott analytic index
under that equivalence.

Let \(\pi:X\to B\) be a proper submersion whose fibers are closed
\(8k\)-dimensional \emph{Spin} manifolds.  Let \(T^{V}X\to X\) be the
vertical tangent bundle, equipped with a Riemannian metric \(g^{T^{V}X}\).
Choose a horizontal distribution \(T^{H}X\subset TX\) and a Riemannian
metric \(g^{TB}\) on \(B\), and set
\[
g^{TX}\ :=\ g^{T^{V}X}\ \oplus\ \pi^{*}g^{TB}.
\]
Let \(\nabla^{TX}\) be the Levi--Civita connection of \(g^{TX}\), and let
\(P:TX\to T^{V}X\) be the orthogonal projection.  Define the natural
connection on \(T^{V}X\) by
\[
\nabla^{T^{V}X}\ :=\ P\circ \nabla^{TX}\circ P .
\]
Assume \(T^{V}X\) is Spin and write \(S^{V}X\to X\) for the vertical spinor
bundle.  The connection \(\nabla^{T^{V}X}\) determines the \(\widehat A\)-form
\(\widehat A(\nabla^{T^{V}X})\in \Omega^{\mathrm{even}}(X)\).

Define the fiber integration with the \(\widehat A\)-correction
\[
\pi_{*}:\ \Omega^{4k-1}(X)\ \longrightarrow\ \Omega^{4k-1}(B),
\qquad
\pi_{*}(\phi)\ :=\ \int_{X/B}\!\widehat A\!\bigl(\nabla^{T^{V}X}\bigr)\wedge \phi .
\]

\paragraph{\emph{Family Dirac operator and index bundle.}}
Fix a real bundle with connection \((E,\nabla^{E})\) over \(X\).  For each
\(b\in B\) let \(F=X_b\) be the fiber, \(\iota_b:F\hookrightarrow X\) the
inclusion, and consider the twisted vertical Spin--Dirac operator
\[
D_b:\ \Gamma\!\bigl(F,\,S^{V}_{+}|_{F}\otimes \iota_b^{*}E\bigr)\ \longrightarrow\
      \Gamma\!\bigl(F,\,S^{V}_{-}|_{F}\otimes \iota_b^{*}E\bigr).
\]
If \(\dim\ker D_b\) and \(\dim\operatorname{coker}D_b\) are locally constant in \(b\),
these spaces form vector bundles over \(B\), and the \emph{index bundle} is
\[
\operatorname{Ind}(D^{V}_{E})\ :=\ [\ker D^{V}_{E}]-[\operatorname{coker}D^{V}_{E}]
\ \in\ KO^{0}(B).
\]
(We use Bott periodicity; the fiber dimension is \(8k\), so the degree shift
is a multiple of \(8\) and the pushforward lands in \(KO^{0}(B)\).)
This gives the analytic pushforward on topological \(KO\):
\[
\pi^{\mathrm{an}}_{!}:\ KO^{0}(X)\longrightarrow KO^{0}(B),
\qquad
\pi^{\mathrm{an}}_{!}(E)\ =\ \operatorname{Ind}(D^{V}_{E}) .
\]

\paragraph{\emph{Connections on the index bundle.}}
Let \(H\) be the (Hilbert) bundle of smooth sections of \(S^{V}X\otimes E\).
Let \(P:H\to \ker D^{V}_{E}\) be the fiberwise orthogonal projection.
Let \(A_{t}\) denote Bismut’s superconnection (built from the horizontal
distribution and \(\nabla^{E}\)), and \(P_{t}\) the spectral projection onto
\(\ker D^{V}_{E}\) defined via \(A_{t}\).
The projected connection on \(\ker D^{V}_{E}\) is denoted \(\nabla^{\mathrm{Ind}}\).

\paragraph{\emph{Bismut--Cheeger transgression form.}}
For the family \(D^{V}_{E}\), there is a canonical odd form
\(\widetilde\eta(D^{V}_{E})\in \Omega^{4k-1}(B)\) such that
\begin{equation}\label{eq:BC-transgression}
d\,\widetilde\eta(D^{V}_{E})
\;=\;
\int_{X/B}\!\widehat A\!\bigl(\nabla^{T^{V}X}\bigr)\wedge \mathrm{ph}\!\bigl(\nabla^{E}\bigr)
\;-\;
\mathrm{ph}\!\bigl(\nabla^{\mathrm{Ind}}\bigr).
\end{equation}

\paragraph{Theorem 7.1}
In the Freed--Lott model of differential \(KO\), the family index theorem on triple
\(\mathcal E=(E,\nabla^{E},\phi)\) with \(\phi\in \Omega^{4k-1}(X)\).
is 
\[
\operatorname{ind}^{\mathrm{an}}_{KO}(\mathcal E)
\;:=\;
\Bigl(\operatorname{Ind}(D^{V}_{E}),\ \nabla^{\mathrm{Ind}},\
      \pi_{*}(\phi)+\widetilde\eta(D^{V}_{E})\Bigr)
\ \in\ \widehat{KO}_{\mathrm{FL}}(B).
\]

\paragraph{Definition 7.2 }
\[
\begin{tikzcd}[column sep=large,row sep=large]
\widehat{KO}_{\mathrm{}}(X) \arrow[r,"f"] \arrow[d,"\mathrm{ind}^{\mathrm{an}}_{\mathrm{}}"'] &
\widehat{KO}_{\mathrm{FL}}(X) \arrow[d,"\mathrm{ind}^{\mathrm{an}}_{\mathrm{FL}}"] \\
\widehat{KO}_{\mathrm{}}(B) &
\widehat{KO}_{\mathrm{FL}}(B) \arrow[l,"g"']
\end{tikzcd}
\]
 define
\[
\mathrm{ind}^{\mathrm{an}}_{\mathrm{}}
\;:=\;
g\ \circ\ \mathrm{ind}^{\mathrm{an}}_{\mathrm{FL}}\ \circ\ f
\;:\;
\widehat{KO}_{\mathrm{}}(X)\longrightarrow \widehat{KO}_{\mathrm{}}(B).
\]

\paragraph{Theorem 7.3.}
In the structured–bundle model, the analytic index
\[
\mathrm{ind}^{\mathrm{an}}_{\mathrm{}}:\ \widehat{KO}_{\mathrm{}}(X)\longrightarrow \widehat{KO}_{\mathrm{}}(B)
\]
is
\begin{equation}\label{eq:ss-an-index}
\mathrm{ind}^{\mathrm{an}}_{\mathrm{}}(\mathcal{E})
=  
\bigl[\operatorname{Ind}(D^{V}_{E}),[\nabla^{\mathrm{Ind}}]\bigr]
\;+\;
\bigl[H,[\nabla^{H}]\bigr]
\;-\;
\bigl[\dim(H),[d]\bigr].
\end{equation}
Here
\[
 \widehat{CS}^{-1}\!\bigl(\widetilde{\eta}(D^{V}_{E})\bigr),
\]
i.e., \(H\) is the (stably trivial) structured bundle whose \(\widehat{CS}\)-class
represents the Cheeger–Bismut eta form \(\widetilde{\eta}(D^{V}_{E})\).
\begin{proof}
This follows from proposition 3.2, and theorem 7.1.
\end{proof}

Let \(\pi:X\to B\) be a smooth submersion with closed \emph{Spin} fibers of
dimension \(8k\).
Given an enriched Spin cycle \((M,f)\) over \(B\) (so \(f:M\to B\) and \(M\) is
closed, equipped with a Riemannian metric and Spin structure), form the
fiber product
\[
\widetilde{M}\ :=\ M\times_{B}X
\ =\ \{(m,x)\in M\times X\mid f(m)=\pi(x)\},
\qquad
p:\widetilde{M}\to M,\ \ \tilde f:\widetilde{M}\to X
\]
for the two projections. Endow \(\widetilde{M}\) with the product Spin structure and the
adiabatic (rescaled) metric \(g_t=p^{*}g^{TM}\oplus t^{-2}\,g^{T(X/B)}\) for \(t>0\).

\paragraph{Definition 7.4}
With \(\pi:X\to B\) a smooth submersion with closed Spin fibers of dimension \(8k\),
and for an enriched Spin cycle \((M,f)\) over \(B\), form
\[
\widetilde{M}=M\times_{B}X,\qquad
p:\widetilde{M}\to M,\ \ \tilde f:\widetilde{M}\to X,
\]
and equip \(\widetilde{M}\) with the adiabatic metric \(g_t=p^{*}g^{TM}\oplus t^{-2}g^{T(X/B)}\).

For a differential \(KO\)-character \(\widehat\kappa\) on \(X\) associated to \((E,\nabla)\),
define \((\pi_{!}\widehat\kappa)(M,f)\) by
\[
(\pi_{!}\widehat\kappa)(M,f)\ :=\
\begin{cases}
\displaystyle \bar\eta\!\left(D^{\,\tilde f^{*}(E,\nabla)}_{\widetilde{M},\,t}\right)\Big|_{t\to 0}
\quad & \text{if } \dim M\equiv 7 \ (\mathrm{mod}\ 8),\\[1.1ex]
\displaystyle \tfrac{1}{2}\,\bar\eta\!\left(D^{\,\tilde f^{*}(E,\nabla)}_{\widetilde{M},\,t}\right)\Big|_{t\to 0}
\quad & \text{if } \dim M\equiv 3 \ (\mathrm{mod}\ 8),\\[1.1ex]
\displaystyle \pi_{!}^{\mathrm{AS}}\!\big(\tilde f^{*}E\big)\in KO^{-(\dim \widetilde{M})}(\mathrm{pt})\cong \mathbb{Z}_{2}
\quad & \text{if } \dim M\equiv 1,2 \ (\mathrm{mod}\ 8),
\end{cases}
\]
where \(D_{\widetilde{M},t}\) is the Spin–Dirac operator on \((\widetilde{M},g_t)\),
\(\bar\eta(D)=\tfrac{\eta(D)+h(D)}{2}\ (\mathrm{mod}\ \mathbb{Z})\), and
\(\pi_{!}^{\mathrm{AS}}\) denotes the Atiyah–Singer \(KO\)-pushforward to a point
(the \(\mathbb{Z}_{2}\)-valued index).

\paragraph{Proposition 7.5}
The pushforward  \(\pi_!:\widehat{KOCH}(X)\to\widehat{KOCH}(B)\) is well defined.
(the following diagram is commutative)
\[
\int_{X/B}\widehat A\big(T(X/B)\big)\wedge(-)\;:\;
\wedge_{\mathrm{integrality}}(X)\longrightarrow \wedge_{\mathrm{integrality}}(B)
\]

\[
\begin{tikzcd}
\widehat{KOCH}(X)
  \arrow[r, "\delta_1"]
  \arrow[d, "\pi_!"']
&
\wedge_{\mathrm{integrality}}(X)
  \arrow[d,     "\displaystyle \int_{X/B}\widehat A(T(X/B))\wedge(-)"]
\\
\widehat{KOCH}(B)
  \arrow[r, "\delta_1"']
&
\wedge_{\mathrm{integrality}}(B)
\end{tikzcd}
\]

we want to show 
The curvature of the pushforward $\pi_!\widehat\kappa$ is
\[
\delta_1(\pi_!\widehat\kappa)
=\int_{X/B}\widehat A\big(T(X/B)\big)\wedge\operatorname{ph}(E,\nabla).
\]

\begin{proof}
By the APS index theorem, for any Spin cobordism $(W,F)$ over $B$ with
pullback $\widetilde W=W\times_B X$ and $\widetilde E=\widetilde F^{*}E$,
\[
(\pi_!\widehat\kappa)(\partial W,F|_{\partial W})
\equiv
\int_{\widetilde W}\widehat A(T\widetilde W)\wedge\operatorname{ph}(\widetilde E)
\quad(\mathrm{mod}\ \mathbb Z).
\]
Let $I$ denote integration along the fibres $\hat\pi:\widetilde W\to W$.
Using $T\widetilde W\cong\hat\pi^{*}TW\oplus F^{*}T(X/B)$ and
$\widetilde E=\widetilde F^{*}E$, we compute
Hence
\[
\begin{aligned}
I\bigl(\widehat A(T\widetilde W)\wedge\operatorname{ph}(\widetilde E)\bigr)
&= I\Bigl(\hat\pi^{*}\widehat A(TW)\wedge
          F^{*}\bigl(\widehat A(T(X/B))\wedge\operatorname{ph}(E,\nabla)\bigr)\Bigr)\\
&= \widehat A(TW)\wedge
   F^{*}I\bigl(\widehat A(T(X/B))\wedge\operatorname{ph}(E,\nabla)\bigr)\\
&= \widehat A(TW)\wedge
   F^{*}\Bigl(\int_{X/B}\widehat A(T(X/B))\wedge\operatorname{ph}(E,\nabla)\Bigr).
\end{aligned}
\]

\[
(\pi_!\widehat\kappa)(\partial W,F|_{\partial W})
\equiv
\int_{\widetilde W}\widehat A(T\widetilde W)\wedge\operatorname{ph}(\widetilde E)
=
\int_W \widehat A(TW)\wedge F^{*}\delta_1(\pi_!\widehat\kappa)
\quad (\mathrm{mod}\ \mathbb Z),
\]
so
\[
\delta_1(\pi_!\widehat\kappa)
=
\int_{X/B}\widehat A(T(X/B))\wedge\operatorname{ph}(E,\nabla),
\]
as claimed.
\end{proof}

\medskip
We have two versions of the family index theorem for Spin fibrations—one in the differential KO-character model and one in the structured-bundle model. By the Cheeger–Bismut adiabatic limit theorem, their pushforwards agree . Moreover, if one could establish the commutativity of the two pushforwards independently, it would yield an alternative proof of the Cheeger–Bismut adiabatic limit theorem.

\paragraph{Theorem 7.6 }
Let \(\pi:X\to B\) be a proper submersion with closed spin fibers of dimension \(8k\).
Then the square
\[
\begin{tikzcd}
\widehat{KO}(X) \arrow[d, "\check{\pi}_!^{\,\mathrm{an}}"'] \arrow[r, "\cong"] &
\widehat{KOCH}(X) \arrow[d, "\pi_!"] \\
\widehat{KO}(B) \arrow[r, "\cong"] &
\widehat{KOCH}(B)
\end{tikzcd}
\]
is commutative.

\begin{proof}

This diagram can be read in two  ways:
\begin{enumerate}[label=\textbf{(\arabic*)}]
\item Apply the isomorphism \(\widehat{KO}(X)\xrightarrow{\ \cong\ }\widehat{KOCH}(X)\) to form the associated differential \(KO\)-characters (see~(5.5)), then push forward on \(KO\)-characters.
\item First pushforward the structured bundle \(\widehat{KO}(X)\xrightarrow{\ \check{\pi}_!^{\,\mathrm{an}}\ }\widehat{KO}(B)\), and then take its associated \(KO\)-character via \(\widehat{KO}(B)\xrightarrow{\ \cong\ }\widehat{KOCH}(B)\).
\end{enumerate}

When the chosen spin cycles have dimension \(4k-1\) , proceeding via \textbf{(1)} yields
$\lim_{t\to 0}\,\bar{\eta}\!\left(D^{\,\tilde{f}^{*}(E,\nabla)}_{\widetilde{M},t}\right)\;(\mathrm{mod}\ \mathbb{Z}) $

Going the other way \textbf{(2)}, we form the \(KO\)-character associated to
\[
\bigl[\operatorname{Ind}(D^{V}_{E}),[\nabla^{\mathrm{Ind}}]\bigr]
\;+\;
\bigl[H,[\nabla^{H}]\bigr]
\;-\;
\bigl[\dim(H),[d]\bigr].
\]
by corollary (6.3)
Evaluating on \(M\) yields
\[
\int_{M}\!\widehat{A}(TM)\wedge f^{*}\widetilde{\eta}\!\left(D^{V}_{E}\right)
\;+\;
\widetilde{\eta}\!\left(D_{M}^{\,f^{*}\!\nabla^{\operatorname{Ind}}}\right)
\;\;(\mathrm{mod}\ \mathbb{Z}) .
\]

Then the equality directly from the Bismut--Cheeger adiabatic limit theorem. in [BC89]
\end{proof}

\section{Examples}
In the final section of this paper, We describe how our results relate to earlier classical developments in the 70s, including Adams’s \(e\)-invariant on framed \((4k-1)\)-manifolds, the relationship to the Cheeger–Simons Pontryagin classes and associated invariants on \(3\)-manifolds, and cobordism invariants for flat bundles discovered by Atiyah–Patodi–Singer.

\subsection{Adam's e invariant}

By the Poincaré polynomial in Anderson–Brown–Peterson. Their splitting result shows that \(\mathrm{MSpin}\) decomposes into a connective \(KO\)-part together with suspensions of Eilenberg–MacLane spectra. In particular, in degrees \(4k-1\) the Hurewicz map
\[
\Omega^{\mathrm{fr}}_{*} \longrightarrow \Omega^{\mathrm{spin}}_{*}
\]
has zero image. It should be noted that in [25], Atiyah mistakenly stated that \(\Omega^{\mathrm{spin}}_{4k-1}=0\). In fact, as computed by Bunke and Naumann in [26], one has \(\Omega^{\mathrm{spin}}_{4k-1}=0\) for \(1\le k\le 9\), but \(\Omega^{\mathrm{spin}}_{39}\cong \mathbb{Z}_{2}\oplus \mathbb{Z}_{2}\neq 0\). However, this does not affect Atiyah’s argument, because the conclusion used there is only that a framed bordism class in degree \(4k-1\) bounds a Spin manifold of dimension \(4k\), rather than that a Spin manifold of dimension \(4k-1\) itself bounds.

Suppose \(E\) of dimension r is a real, flat vector bundle with connection \(\nabla^{0}\), so its Pontryagin character appears only in degree \(0\) (recording the rank). To apply the \(KO\)-character associated with \((E,\nabla^{0})\) to an enriched Spin cycle, the degree may be \(8k+1\), \(8k+2\), \(8k+3\), or \(8k+7\). We consider the case in which this \(KO\)-character is evaluated on an \(8k+3\) or \(8k+7\) enriched Spin cycle \(M\) arising from the Hurewicz map; that is, \(M\) is a framed bordism cycle  carrying a preferred Spin structure induced by a trivialization of the stable tangent bundle. The enriched connection on \(M\) is induced by this trivialization as well. Since \(M\) bounds a manifold \(W\), the value of the \(KO\)-character on \(M\) is equal to the evaluation on \(W\) of \(\widehat{A}(W)\) multiplied by the Pontryagin character of \(f^{*}E\) with \(f^{*}\nabla^{0}\), yielding a constant times the characteristic number \(\widehat{A}(W)\) in \(\mathbb{Q}/\mathbb{Z}\). If \(\dim M = 8k+3\), this value is further divided by \(2\).

$\pi$ denotes the choice of framing. 

\[
\hat{\kappa}_{(E,\nabla^{0})}(M,\pi)=
\begin{cases}
r\widehat{A}(TN)[N], & \text{if }  \text{ n=8k+7},\\[4pt]
\dfrac{r}{2}\,\widehat{A}(TN)[N], & \text{if }  \text{ n=8k+3}.
\end{cases}
\]

By the Pontryagin–Thom construction, the isomorphism between  stable homotopy groups of spheres and framed cobordism. In particular, a parallelism \(\pi\) on a manifold \(M\) induces a framing on \(M\) and hence defines an element
\[
[M,\pi]\in \pi^{S}_{n}, \qquad \text{where } n=\dim M.
\]
The Adams \(e\)-invariant for \(n=4k-1\) is a homomorphism
\[
e \colon \pi^{S}_{4k-1} \longrightarrow \mathbb{Q}/\mathbb{Z}.
\] ( we give a quick recount of the construction of this homomorphism)

 Let \(m,n \ge 1\), and let \(f \colon S^{2n+2m-1} \to S^{2n}\) be a pointed map. Let
\[
X = X_f = S^{2n} \cup_f e^{2n+2m}
\]
be the mapping cone of \(f\), let \(i \colon S^{2n} \hookrightarrow X\) be the inclusion, and let
\[
\pi \colon X \longrightarrow X / S^{2n} \cong S^{2n+2m}.
\]
Then the cofibration
\[
S^{2n+2m-1} \xrightarrow{\,f\,} S^{2n} \xrightarrow{\,i\,} S^{2n} \cup_f e^{2n+2m} \xrightarrow{\,\pi\,} S^{2n+2m}
\]
 induces a long exact sequence in reduced \(K\)-theory. Since the \(K\)-theory of spheres is concentrated in even degrees, the \(K\)-theory degree of \(f\), i.e., \(\widetilde{K}(f)\), is zero. 
Since \(\widetilde{K}(f)=0\), we obtain a short exact sequence
\begin{equation}\label{eq:short-exact}
0 \longrightarrow \widetilde{K}\!\left(S^{2n+2m}\right)
\xrightarrow{\ \pi^{*}\ } \widetilde{K}\!\left(S^{2n} \cup_f e^{2n+2m}\right)
\xrightarrow{\ i^{*}\ } \widetilde{K}\!\left(S^{2n}\right)
\longrightarrow 0 .
\end{equation}

\medskip

Let \(i_{2n}\) be a generator of \(\widetilde{K}(S^{2n})\) and \(i_{2n+2m}\) a generator of \(\widetilde{K}(S^{2n+2m})\). Choose an element
\[
a \in \widetilde{K}\!\left(S^{2n} \cup_f e^{2n+2m}\right)
\quad\text{such that}\quad
i^{*}(a)=i_{2n},
\]
and let \(b=\pi^{*}(i_{2n+2m}) \in \widetilde{K}\!\left(S^{2n} \cup_f e^{2n+2m}\right)\).
Then, for any \(k\), we have the Adam's operation
\[
\psi^{k}(a)=k^{n}\cdot a+\mu_{k}\cdot b.
\]
Since the Adams operations commute, we must have \(\psi^{\ell}(\psi^{k}(a))=\psi^{k}(\psi^{\ell}(a))\), and hence
\[
k^{n}\!\left(k^{m}-1\right)\mu_{\ell}=\ell^{\,n}\!\left(\ell^{\,m}-1\right)\mu_{k}
\quad\text{for any } k,\ell.
\]
This shows that the rational number
\[
e(f):=\frac{\mu_{k}}{k^{n}\!\left(k^{m}-1\right)} \in \mathbb{Q}
\]
is independent of \(k\). It might, however, depend on our choice of \(a\). If we change \(a\) by a multiple of \(b\), then \(e(f)\) changes by an integer. (For \(a' = a + p\cdot b\), we obtain \(e'(f)=e(f)+p\).) Thus \(e(f)\) is well defined as an element of \(\mathbb{Q}/\mathbb{Z}\).
Thus an assignment
\[
\bigl(f \colon S^{2n+2m-1} \to S^{2n}\bigr) \longmapsto e(f) \in \mathbb{Q}/\mathbb{Z}.
\]

\begin{theorem}[Atiyah–Patodi–Singer]
The \(e\)-invariant of a framed \(4k-1\)-dimensional manifold \((M,\pi)\) is given by
\[
e(M,\pi)
=\varepsilon(k)\left(
\frac{\eta(D_M)+h(D_M)}{2}
\;+\;
\int_{M}\widetilde{\widehat{A}}\!\bigl(TM,\nabla^{TM},\nabla^{\pi}\bigr)
\right)\in \mathbb{Q}/\mathbb{Z},
\]
where \(\varepsilon(k)=1\) if \(k\) is even and \(\varepsilon(k)=\tfrac{1}{2}\) if \(k\) is odd.
\end{theorem}

This implies that, in the special case of applying the \(KO\)-character associated with a flat bundle to an enriched Spin cycle arising from the Hurewicz map from framed bordism, the invariant is independent of the geometric data as we have a topological description thanks to APS and Adam.

\subsection{invariants on the 3 spin manifold }

Let \((E,\nabla^{E})\) be a real vector bundle with connection over \(X\), and let
\(f\colon M^{3}\to X\) be a Spin bordism cycle. For two auxiliary connections
\(\nabla^{0},\nabla^{1}\) on \(TM\), the variation of the associated differential \(KO\)-character satisfies

Since the Spin cobordism group in degree \(3\) vanishes, i.e., \(\Omega^{\mathrm{spin}}_{3}=0\).

Since \(\widehat{A}=1-\tfrac{1}{24}p_{1}\) in degree \(4\), it follows that
\[
\widetilde{\widehat{A}}\!\bigl(TM;\nabla^{0},\nabla^{1}\bigr)
= -\frac{1}{24}\,\widetilde{p}_{1}\!\bigl(TM;\nabla^{0},\nabla^{1}\bigr),
\]
and hence
\[
\hat{\kappa}_{(E,\nabla^{E})}(M,\nabla^{0})
-
\hat{\kappa}_{(E,\nabla^{E})}(M,\nabla^{1})
=
-\frac{1}{24}\int_{M}
\widetilde{p}_{1}\!\bigl(TM;\nabla^{0},\nabla^{1}\bigr)\wedge
\operatorname{ph}\!\bigl(f^{*}\nabla^{E}\bigr).
\]

\(\widehat{p}_{1}(E,\nabla^{E,1})\) denotes the canonical Cheeger–Simons lift of the first Pontryagin class, and \(\widehat{\operatorname{ph}}\!\big(f^{*}E,\,f^{*}\nabla^{E}\big)\) denotes the Cheeger–Simons lift of the Pontryagin character form associated to \(f^{*}\nabla^{E}\).

\noindent
Here \(\widetilde{p}_{1}\) denotes the Chern–Simons class associated with the first Pontryagin class. The variation formula for the associated Cheeger–Simons character \(\widehat{p}_{1}\) of a general vector bundle is
\[
\bigl(\,\widehat{p}_{1}(E,\nabla^{E,1})-\widehat{p}_{1}(E,\nabla^{E,0})\,\bigr)[M]
=\int_{M}\widetilde{p}_{1}\!\bigl(E,\nabla^{E,0},\nabla^{E,1}\bigr).
\]

Consider the natural product structure on differential cohomology. Then

\[
\widehat{H}^{p}(M)\times \widehat{H}^{q}(M)\longrightarrow \widehat{H}^{p+q}(M).
\]

\[
\hat{\kappa}_{(E,\nabla^{E})}(M,\nabla^{0})
-
\hat{\kappa}_{(E,\nabla^{E})}(M,\nabla^{1})
=
-\frac{1}{24}\,
\big\langle
\big(\widehat{p}_{1}(TM,\nabla^{1})-\widehat{p}_{1}(TM,\nabla^{0})\big)
\,\widehat{\cup}\,
\widehat{\operatorname{ph}}\!\big(f^{*}E,f^{*}\nabla^{E}\big),
\,[M]\big\rangle
\]
\[
\in \mathbb{R}/\mathbb{Z}.
\]

\subsection{invariants on the flat bundle}
Since unitary representations of the fundamental group are in one-to-one correspondence with flat Hermitian bundles [27, p.55], we may represent a flat Hermitian bundle by its holonomy representation
\[
\alpha:\pi_{1}(M)\longrightarrow U(n).
\]

Let \(\alpha:\pi_{1}(M)\to U(n)\) be a unitary representation of the fundamental group. Then
\[
F_{\alpha}=\widetilde{M}\times_{\alpha}\mathbb{C}^{n}\longrightarrow M
\]
is a flat Hermitian vector bundle with holonomy \(\alpha\).

\(D_{\alpha}\) denotes the Dirac operator twisted by the flat Hermitian bundle.

\begin{theorem}[APS Cobordism invariants on Flat bundle ]

\end{theorem}
\bigskip

The \(\xi\)-invariant
\[
\xi_{\alpha}(D)=\frac{\eta+h}{2}(D_{\alpha})-n\cdot\frac{\eta+h}{2}(D)\;\in\;\mathbb{R}/\mathbb{Z}
\]
is a cobordism invariant in the sense that \(\xi_{\alpha}(D)=0\) if there exists a compact manifold \(N\) with \(M=\partial N\) such that \(D\) extends to an operator on \(N\), and \(\alpha\) extends to a representation of \(\pi_{1}(N)\).

The representation \(\alpha\) defines a class \([\alpha]\in K^{-1}(M;\mathbb{R}/\mathbb{Z})\), and the symbol of \(D\) gives \(\sigma\in K^{1}(TM)\).

There is a natural topological product 
\[
K^{-1}(M;\mathbb{R}/\mathbb{Z}) \otimes K^{1}_{c}(TM)
\;\xrightarrow{\ \smile\ }\;
K^{0}_{c}(TM;\mathbb{R}/\mathbb{Z}).
\]

Then there exists a topological index \(\operatorname{Ind}_{[\alpha]}(\sigma)\in\mathbb{R}/\mathbb{Z}\), and
\[
\xi_{\alpha}(D_{M})=\operatorname{Ind}_{[\alpha]}(\sigma).
\]

We now explain the connection between differential K character induced by \(K^{-1}(X;\mathbb{R}/\mathbb{Z})\)  and  this cobordism invariant from APS.

Recall we have the natural topological pairing through slant product.
\[
K^{-1}(X;\mathbb{R}/\mathbb{Z})\otimes KU_{\mathrm{odd}}(X)\longrightarrow \mathbb{R}/\mathbb{Z}
\]
This pairing admits an analytic formula in the following way 

An element of \(KU_{\mathrm{odd}}(X)\) can be represented by the Conner–Floyd class of a complex bordism cycle \((M,f)\).

Any element of \(K^{-1}(X;\mathbb{R}/\mathbb{Z})\) is represented by 
\[
[E,\nabla]- (r) ,
\] 
where \((E,\nabla)\) is  some rank r hermitian vector bundle with vanishing chern-weil forms, \((r)\) denotes the trivial bundle \(\mathbb{C}^{r}\) with the trivial connection. The pairing, by differential K Character in [SS18] is 

\[
\big\langle [E,\nabla] - (r),\ (M,f)\big\rangle
\;\equiv\;
\bar\eta\!\left(D^{\,\nabla}_{M,f^*E}\right)
\;-\; r\,\bar\eta(D_M)
\]

Our pairing can be viewed as a generalization of the APS cobordism invariant. In case of thereom 8.2, here \(D\) is the \(\mathrm{Spin}^{c}\) Dirac operator, and \(\alpha\) denotes its twist by the pullback \(f^{*}(E,\nabla^{E})\) of the Hermitian bundle with vanishing chern-weil forms \((E,\nabla^{E})\) along \(f\colon M\to X\).

It should be noted that our construction admits bundles with connections whose Chern–Weil characteristic forms vanish (i.e., they lie in the kernel of the curvature map in the differential K theory hexagon), without requiring the connection to be flat. For example, identify \(S^{3}\simeq SU(2)\), take the trivial rank-\(2\) bundle \(E=S^{3}\times\mathbb{C}^{2}\) with structure group \(SU(2)\subset U(2)\), and let \(\theta=g^{-1}dg\) be the left-invariant Maurer–Cartan \(\mathfrak{su}(2)\)-valued \(1\)-form, so \(d\theta+\theta\wedge\theta=0\). For any constant \(c\in\mathbb{R}\setminus\{0,1\}\), define
\[
A:=c\,\theta,\qquad \nabla:=d+A,\qquad
F_{\nabla}=dA+A\wedge A=(c^{2}-c)\,\theta\wedge\theta\neq 0.
\]
Since \(A,F_{\nabla}\in\mathfrak{su}(2)\) are traceless and \(\dim S^{3}=3\),
\[
\operatorname{ch}_{1}(E,\nabla)=\frac{i}{2\pi}\operatorname{Tr}(F_{\nabla})=0,
\qquad
\operatorname{ch}_{k}(E,\nabla)=0\;\;(k\ge2),
\]
hence
\[
\operatorname{ch}(E,\nabla)=\mathrm{rk}(E)=2,
\]
even though \(F_{\nabla}\ne 0\).  This is an example of a non-flat bundle with vanishing chern-weil forms.

Our result further suggests that, after restricting to flat bundle in the \(K^{-1}(X;\mathbb{R}/\mathbb{Z})\)  and taking \(D\) to be the \(\mathrm{Spin}^{c}\) Dirac operator, the two topological pairings agree
\[
\langle\,\cdot,\cdot\,\rangle_{}\;:\;
K^{-1}(X;\mathbb{R}/\mathbb{Z})\otimes KU_{\mathrm{odd}}(X)\longrightarrow \mathbb{R}/\mathbb{Z},
\qquad
\langle\,\cdot,\cdot\,\rangle_{}\;:\;
K^{-1}(X;\mathbb{R}/\mathbb{Z})\otimes K^{1}_{c}(TM)\longrightarrow \mathbb{R}/\mathbb{Z}
\]

\appendix
\section{From ordinary characters to Characteristic Variety theorem on ordinary Cohomology}

In this section, we present a Characterisitc Variety theorem like method for lifting a differential character to a cocycle in ordinary cohomology. From the introduction, we see that Morgan–Sullivan showed that each class in ordinary cohomology is isomorphic to a map from homology with $\mathbb{Q}/\mathbb{Z}$ coefficients to $\mathbb{Q}/\mathbb{Z}$ that is continuous with respect to the profinite topology (and thus can be lifted to homology with $\mathbb{Q}$ coefficients). In terms of the characteristic variety theorem, this means that a cohomology class is completely determined by its $\mathbb{Z}$ and $\mathbb{Z}/k$-periods.

We now begin our construction. \[
\text{We start with } f \text{ as a differential character, } f \in \widehat{H}^{\,k}(X).
\]Since $\mathbb{R}$ is divisible, we can find $T \in C^{k-1}(M;\mathbb{R})$ with $\widetilde{T}\!\mid_{Z_{k-1}(M)} = f$. It is easily seen that $\omega_f - \delta T \in C^k(M;\mathbb{Z})$ and that it is closed. We define the following map on $\mathbb{Z}_k$-manifolds :
\[
\operatorname{Inv}_{\mathbb{Z}_k}\colon H_{i}\!\left(X;\,\mathbb{Z}/k\mathbb{Z}\right)\longrightarrow \mathbb{Z}/k\mathbb{Z}.
\]

\[
\operatorname{Inv}_{\mathbb{Z}_k}\colon M \longmapsto \frac{1}{k}\int_M \bigl(\delta {T}-\omega_f\bigr) \bmod \mathbb{Z}\
.
\]

We now verify this map is well defined; specifically, $\operatorname{Inv}_{\mathbb{Z}_k}$ evaluates to $0$ on any $\mathbb{Z}_k$-boundary.

Recall that a $\mathbb{Z}_k$-manifold $M$ is called a $\mathbb{Z}_k$-boundary if its boundary consists of $k$ copies of the Bockstein, each of which bounds a manifold $W$, and $M \sqcup kW$ is a boundary.

Since $\operatorname{Inv}_{\mathbb{Z}_k}$ on $M \sqcup kW$ is $0$ by Stokes’ theorem, we have
\[
0 \equiv \operatorname{Inv}_{\mathbb{Z}_k}(M) \;+\; \int_{W} (\delta T - \omega_f) \pmod{\mathbb{Z}}.
\]
But $\delta T - \omega_f$ is an cochain in integer coefficient, so $\displaystyle \int_{W} (\delta T - \omega_f) \in \mathbb{Z}$. Hence the second term vanishes modulo $\mathbb{Z}$, and therefore $\operatorname{Inv}_{\mathbb{Z}_k}(M)=0$ whenever $M$ is a $\mathbb{Z}_k$-boundary.

Taking direct limits over $k$, we obtain
\[
\operatorname{Inv}_{\mathbb{Q}/\mathbb{Z}}
=\varinjlim_{k}\,\operatorname{Inv}_{\mathbb{Z}/k\mathbb{Z}}
\qquad\text{and}\qquad
H_{i}\!\left(X;\,\mathbb{Q}/\mathbb{Z}\right)
\cong \varinjlim_{k}\,H_{i}\!\left(X;\,\mathbb{Z}/k\mathbb{Z}\right).
\]

\begin{center}
\[
\begin{tikzcd}[column sep=4em,row sep=3em,baseline=(current bounding box.center)]
H_{i}(X;\mathbb{Q}) \arrow[r] \arrow[d, "\pi_{*}"'] &
\mathbb{Q} \arrow[d, "\pi"] \\
H_{i}(X;\mathbb{Q}/\mathbb{Z}) \arrow[r, "\operatorname{Inv}_{\mathbb{Q}/\mathbb{Z}}"] &
\mathbb{Q}/\mathbb{Z}
\end{tikzcd}
\]
\end{center}

recall \paragraph*{Morgan–Sullivan (1974) and Pontryagin duality.}
Recall that the Pontryagin dual of a topological abelian group \(G\) is
\(\widehat{G}:=\Hom_{c}(G,S^{1})\) (continuous characters).
Morgan–Sullivan [MS74] prove that for any finite CW–complex \(X\) and each \(i\)
there is a natural isomorphism
\[
\left[
\vcenter{\hbox{%
\begin{tikzcd}[column sep=4em,row sep=3em,baseline=(current bounding box.center)]
H_{i}(X;\mathbb{Q}) \arrow[r] \arrow[d, "\pi_{*}"'] &
\mathbb{Q} \arrow[d, "\pi"] \\
H_{i}(X;\mathbb{Q}/\mathbb{Z}) \arrow[r] &
\mathbb{Q}/\mathbb{Z}
\end{tikzcd}%
}}
\right]
\ \;\cong\;\ 
\Hom_{c}\!\big(H_{i}(X;S^{1}),\,S^{1}\big),
\]
where \(\pi:\mathbb{Q}\twoheadrightarrow \mathbb{Q}/\mathbb{Z}\) is the quotient and
\(\pi_{*}\) is the induced map on homology.  Thus the “diagram–dual’’ on the left
is canonically identified with the Pontryagin dual of the \(S^{1}\)-homology group.

Consequently, in the case where the Anderson–dual theory satisfies
\(Dh_{i}(X;{-}) \cong H_{i}(X;{-})\) (e.g. when we are working with ordinary homology),
the \(\mathbb{Q}/\mathbb{Z}\)-bordism isomorphism for the Anderson dual recovers the
classical Pontryagin duality isomorphism above.

\section{Differential Euler character at odd dimensional cycle}

Recall that the differential \(KO\)-character evaluated on Spin cycles of
dimensions \(8K+1\) and \(8K+2\) is \(\mathbb{Z}_{2}\)-torsion; consequently,
there is no associated characteristic form in these degrees. The same
\(\mathbb{Z}_{2}\)-phenomenon already appears for the differential Euler
character in the Cheeger–Simons setting.

Let \(X\) be a smooth manifold and \(E\to X\) a real, oriented vector bundle of
rank \(m\). Write \(e(E)\in H^{m}(X;\mathbb{Z})\) for the Euler class. If \(m\)
is odd, then
\begin{equation*}
2\,e(E)=0 .
\end{equation*}
Let \(\nabla^{E}\) be a metric connection on \(E\). Denote by
\(e(E,\nabla^{E})\in \Omega^{m}(X)\) its Chern–Weil representative and by
\(\widehat e(E,\nabla^{E})\in \widehat H^{\,m-1}(X;\mathbb{R}/\mathbb{Z})\) its
differential refinement. When \(m\) is odd we have \(e(E,\nabla^{E})=0\), hence
\(\widehat e(E,\nabla^{E})\) is independent of \(\nabla^{E}\).

We give the following theorem without proof.

Let \(\pi\colon S(E)\to X\) be the unit sphere bundle, and let \(TS(E)\) be the
vertical tangent bundle (along the spheres) of rank \(m-1\); denote its Euler
form by \(e(TS(E))\). For any cycle \(z\in Z_{m-1}(X)\) there exist
\(y\in Z_{m-1}(S(E))\) and \(w\in C_{m}(X)\) with
\[
z=\pi_{*}(y)+\partial w .
\]

\paragraph{Theorem}[17]
    If \(m\) is odd, then
\[
\big\langle \widehat e(E,\nabla^{E}),\, z \big\rangle
= -\frac{1}{2}\int_{y} e\!\big(TS(E)\big).
\]
In particular,
\[
2\,\widehat e(E,\nabla^{E})=0 \quad \text{in } \widehat H(X).
\]

\end{document}